\documentclass[10pt]{article}
\usepackage{latexsym, amsfonts, amsmath, amssymb, amscd, epsfig}

\numberwithin{equation}{section}
\newtheorem{thm}{Theorem}[section]
\newtheorem{defn}[thm]{Definition}
\newtheorem{pro}[thm]{Proposition}
\newtheorem{lem}[thm]{Lemma}
\newtheorem{asu}[thm]{Assumption}

\newtheorem{re}{Remark}[section]
\newenvironment{pf}{{\noindent \it \bf Proof:}}{{\hfill$\Box$}\\}
\newcommand{\R}{\mathbb{R}}
\newcommand{\p}{\partial}

\usepackage{mathrsfs}
\usepackage{amsmath}
\usepackage{color}

\allowdisplaybreaks

\begin{document}
\title{Existence and Blow-up of solutions for Stochastic Modified Two-component Camassa-Holm System}
\author{Wujun Lv\\
{\small Department of Statistics, College of Science, Donghua University}\\
{\small 201620, Shanghai, P. R. China}\\
{\small lvwujunjaier@gmail.com}\\
Xing Huang\thanks{Corresponding author}\\
{\small Center for Applied Mathematics, Tianjin University}\\
{\small 300072, Tianjin, P. R. China}\\
{\small xinghuang@tju.edu.cn}\\
}

\date{}
\maketitle
\bigskip \noindent \textbf{Abstract:} In this paper, we consider the modified two-component Camassa-Holm System with multiplicative noise. For these SPDEs, we first establish the local existence and pathwise uniqueness of the pathwise solutions in Sobolev spaces $H^{s}\times H^{s}, s>\frac{3}{2}$. Then we
show that strong enough noise can actually prevent blow-up with probability 1. Finally, we analyse the effects of weak noise and present conditions on the initial data that lead to the global existence and the blow-up
in finite time of the solutions, and their associated probabilities are also obtained. 

\section{Introduction}
Consider the following integrable two-component Camassa-Holm (CH2) shallow water system
\begin{equation*}\label{CH2}
\left\{
\begin{array}{l}
(u-u_{xx})_{t}+3uu_{x}-2u_{x}u_{xx}-uu_{xxx}+\rho\rho_x=0,\ \  t>0,~x\in\R, \\
\rho_t+(\rho u)_x=0,\ \  t>0,~x\in\R, \\
u(0,x)=u_0(x),\ \   x\in\R, \\
\rho(0,x)=\rho_0(x), \ \  x\in\R.
\end{array}
\right.
\end{equation*}
This system appears initially in \cite{PP}. In 2008, it was derived by Constantin and Ivanov in \cite{AR}, which provided a demonstration about its derivation
in view of the fluid shallow water theory from the hydrodynamic point of view. Similar to the Camassa-Holm equation, this system possesses the peakon, multi-kink solutions and the bi-Hamiltonian structure \cite{MSY,G} and is
integrable. Well-posedness and wave breaking mechanism were discussed in \cite{JOZ,GY,ZY} and the existence of global solutions was analyzed in \cite{AR,GY,Z2}.

Obviously, under the constraint of $\rho(x, t)=0$, this system reduces to the cerebrated Camassa-Holm (CH) equation, which was derived physically by Camassa and Holm in \cite{RD} (found earlier by Fokas and
Fuchssteiner \cite{BA} as a bi-Hamiltonian generalization of the KdV equation) by approximating directly the Hamiltonian for Euler's equation in the shallow water region with $u(x, t)$ representing the free surface above a flat bottom. CH equation is completely integrable \cite{AH,AVR} and has infinitely many conservation laws \cite{J}. Local well-posedness for the initial data $u_{0}\in H^{s}$ with $s>3/2$ was proved in \cite{AJ, YP}. One of the remarkable features of the CH equation is the presence of breaking waves as well as global solutions in time. Wave breaking for a large class of initial data has been established in \cite{AJ,YP,AJ1,HP,Z,Z1}. Global solutions were also explored in \cite{AJ,AJ1}. The solitary waves of the CH equation are peaked solitons and are orbitally stable \cite{AW}. If $\rho(x,t)\neq0$, this CH2 system is actually an extension of the CH equation.


However, a modified version of the two-component Camassa-Holm (MCH2) system allows a dependence on the average density $\overline{\rho}$ as well as the pointwise density $\rho$, and it is written as
\begin{equation}\label{MCH2}
\left\{
\begin{array}{l}
(u-u_{xx})_{t}+3uu_{x}-2u_{x}u_{xx}-uu_{xxx}+\rho\overline{\rho}_x=0,\ \  t>0,~x\in\R, \\
\rho_t+(\rho u)_x=0,\ \  t>0,~x\in\R, \\
u(0,x)=u_0(x),\ \   x\in\R, \\
\rho(0,x)=\rho_0(x), \ \  x\in\R,
\end{array}
\right.
\end{equation}
where $u$ denotes the velocity field, $\rho=(1-\partial_{x}^{2})(\overline{\rho}-\overline{\rho}_{0})$ with some constant $\overline{\rho}_{0}$. This system was introduced by Holm et al. in \cite{DLC}, and it does admit peaked solutions in the velocity and average density. Many authors analytically identified the steepening mechanism that allows the singular solutions to emerge from smooth spatially confined initial data. They found that wave breaking in the fluid velocity does not imply singularity in the pointwise density $\rho$ at the point of vertical slope. \eqref{MCH2} may not be integrable unlike the CH2 system. The characteristic is that it will amount to strengthening the norm for $\overline{\rho}$ from $L^{2}$ to $H^{1}$ in the potential energy term. Letting $\gamma=\overline{\rho}-\overline{\rho}_{0}$, it leads to the conserved quantity $\int_{\mathbb R}\|u\|_{H^{1}}^{2}+\|\gamma\|_{H^{1}}^{2}dx$, which is absent in the CH2 system. This property inspired a series of interesting works for a deep insight into the MCH2 system in the recent years.
The Cauchy problem of \eqref{MCH2} has been studied in many works \cite{CKZ,ZM,DLC,JZ}. It has been shown that this system is locally well-posed on the line \cite{CKZ} and on the circle \cite{JZ}. Moreover, the authors presented several blow-up results \cite{CKZ,ZM,JZ,ZML}. In addition, basing on a conserved quantity, the authors established the  global existence results for strong solutions to the system \cite{TY}.


Before introducing our model, we recall some theory of infinite dimensional stochastic analysis.
Let
\begin{equation*}\label{S}
\mathcal{S}=(\Omega,\mathcal{F},\mathbb{P}, \{\mathcal{F}_{t}\}_{t\geq0},\mathcal{W}_{1}, \mathcal{W}_{2}),
\end{equation*}
where $(\Omega,\mathcal{F},\mathbb{P}, \{\mathcal{F}_{t}\}_{t\geq0})$ is a complete filtration probability space,
and $\mathcal{W}_{1}, \mathcal{W}_{2}$ are two cylindrical Wiener process on some separable Hilbert space $U$ and $d\langle \mathcal{W}_{1}, \mathcal{W}_{2}\rangle_{t}=\kappa dt, -1\leq\kappa\leq1$.
To be precise, we consider a separable Hilbert space $U$ as well as a larger Hilbert space $U_{0}$ such that the canonical embedding $U\hookrightarrow U_{0}$ is Hilbert-Schmidt. Therefore we have
\begin{equation*}
\mathcal{W}_{i}=\sum_{k=1}^{\infty}W_{k}^{i}e_{k}\in C([0, \infty), U_{0}),~~i=1,2,
\end{equation*}
where $\{W_{k}^{i}\}_{k\geq1}$ is a sequence of mutually independent one-dimensional Brownian motions and $\{e_{k}\}_{k\in\mathbb N}$
is a complete orthonormal basis of $U$.

To define the It\^{o} stochastic integral
\begin{equation*}\label{W}
\int_{0}^{t}Gd\mathcal{W}_{i}=\sum_{k=1}^{\infty}\int_{0}^{t}Ge_{k}dW_{k}^{i}, ~~i=1,2
\end{equation*}
on $H^{s}$, it is required in \cite{DG,PR} for the predictable stochastic process $G$ to take values in the space of $L_{2}(U; H^{s})$, the Hilbert-Schmidt operators from $U$ to $H^{s}$. We have
\begin{equation*}
\bigg(\int_{0}^{t}Gd\mathcal{W}_{i}, v\bigg)_{H^{s}}=\sum_{k=1}^{\infty}\int_{0}^{t}(Ge_{k}, v)_{H^{s}}dW_{k}^{i},~~i=1,2.
\end{equation*}
Moreover, the Burkholder-Davis-Gundy inequality
\begin{equation*}\label{WBDG}
\mathbb E\bigg(\sup_{t\in[0, T]}\left\|\int_{0}^{t}Gd\mathcal{W}_{i}\right\|_{H^{s}}^{p}\bigg)\leq C(p,s)\mathbb E\bigg(\int_{0}^{T}\|G\|^{2}_{L_{2}(U; H^{s})}ds\bigg)^{\frac{p}{2}}, ~p\geq1,~ i=1,2
\end{equation*}
holds for some constant $C(p,s)>0$.

In this paper, we are interested in stochastic variants of the MCH2 system to model energy consuming/exchanging mechanisms in \eqref{MCH2} that are driven by external stochastic influences. Adding multiplicative
noise has also been connected to the prevailing hypotheses that the onset of turbulence in fluid models involves randomness \cite{ZMF,SA,W}.
Precisely, we consider stochastic modified two-component Camassa-Holm (SMCH2) system
\begin{equation}\label{SMCH2}
\left\{
\begin{array}{l}
(u-u_{xx})_{t}+3uu_{x}-2u_{x}u_{xx}-uu_{xxx}+\rho\overline{\rho}_x=(1-\partial_{x}^{2})h_{1}(t,u,\rho)\dot{\mathcal{W}_{1}},\ \  t>0,~x\in\R, \\
\rho_t+(\rho u)_x=(1-\partial_{x}^{2})h_{2}(t,u,\rho)\dot{\mathcal{W}_{2}},\ \  t>0,~x\in\R, \\
u(0,x)=u_0(x),\ \   x\in\R, \\
\rho(0,x)=\rho_0(x), \ \  x\in\R,
\end{array}
\right.
\end{equation}
where  $h_{1}(t,u,\rho), h_{2}(t,u,\rho)$ are typically nonlinear functions.

Let $\gamma=\bar{\rho}-\bar{\rho}_0$, then $(1-\partial_{x}^{2})^{-1}\rho=\gamma$. Notice that the deterministic MCH2 type equations with the weakly dissipative term $\lambda_{2}(1-\partial_{x}^{2})h_{1}(t,u,\rho), \lambda_{2}(1-\partial_{x}^{2})h_{2}(t,u,\rho)$ have been introduced and studied by many scholars \cite{JM,SZ,SZ1}. In order to model more general random energy exchange, we consider the possibly nonlinear noise term $(1-\partial_{x}^{2})h_{1}(t,u,\rho)\dot{\mathcal{W}_{1}}, (1-\partial_{x}^{2})h_{2}(t,u,\rho)\dot{\mathcal{W}_{2}}$ in \eqref{SMCH2}, which will be used to compare with deterministic weakly dissipative MCH2 type equations.
%
%

In \eqref{SMCH2}, the operator $(1-\partial_{x}^{2})^{-1}$ can be expressed by it's associated Green's function $G(x)=e^{-|x|}/2$ with
\begin{equation*}\label{G}
[(1-\partial_{x}^{2})^{-1}f](x)=[G\ast f](x)=\frac{1}{2}\int_{\mathbb R}e^{-|x-y|}f(y)dy.
\end{equation*}
So the system \eqref{SMCH2} is equivalent to tha following one
\begin{eqnarray}\label{SMCH21}
\left\{
\begin{array} {l}
d u+[uu_x+F_{1}(u,\gamma)]d t=h_{1}(t, u,\gamma)d\mathcal{W}_{1},~t>0,~x\in\R, \\
d \gamma+[u\gamma_x+F_{2}(u,\gamma)]d t=h_{2}(t, u,\gamma)d\mathcal{W}_{2},~t>0,~x\in\R, \\
u(0,x)=u_0(x),~x\in\R,\\
\gamma(0,x)=\gamma_0(x),~x\in\R,
\end{array}
\right.
\end{eqnarray}
where $F_{1}(u,\gamma)=\p_x(1-\partial_{x}^{2})^{-1}(u^2+\frac{1}{2}u_x^2+\frac{1}{2}\gamma^2-\frac{1}{2}\gamma_x^2)$, $F_{2}(u,\gamma)=(1-\partial_{x}^{2})^{-1}((u_x\gamma_x)_x+u_x\gamma)$.

The purpose of this paper is as follows:

$\bullet$ The first goal of the present paper is to analyze the existence and uniqueness of pathwise solutions and to determine possible blow-up criterion for the Cauchy problem \eqref{SMCH21}. Under generic assumptions on $h_{1}(t, u,\gamma), h_{2}(t, u,\gamma)$, we will show that \eqref{SMCH21} has a local unique pathwise solution(see Theorem \ref{Maximal solutions} below).

$\bullet$ The second goal of this work is to study the case of strong nonlinear noise and consider its effect. As we will see in \eqref{Blow-up criterion} below, for the solution to \eqref{SMCH21}, its $H^{s}\times H^{s}$-norm blows up if and only if its $W^{1,\infty}\times W^{1,\infty}$-norm blows up. This suggests choosing a noise coefficient involving the $W^{1,\infty}\times W^{1,\infty}$-norm of $(u,\gamma)$. Therefore in this work we consider the case that $h_{1}(t, u, \gamma)d \mathcal{W}_{1}=h_{1}(t, u, \gamma)d (\sum_{i=1}^{\infty}W_{i}(t)e_{i})=a(t)(1+\|u\|_{W^{1,\infty}}+\|\gamma\|_{W^{1,\infty}})^{\theta}ud W_{1}$, where $\{e_{i}\}_{i\in \mathbb N^{*}}$ denote an orthonormal basis, $\{W_{i}\}_{i\in \mathbb N^{*}}$ is a family of independent standard real-valued Wiener processes. Similarly, we consider the case that $h_{2}(t, u, \gamma)d \mathcal{W}_{2}=a(t)(1+\|u\|_{W^{1,\infty}}+\|\gamma\|_{W^{1,\infty}})^{\theta}\gamma d W_{1}$, where $\theta>0, 0<a_{*}\leq a^{2}(t)\leq a^{*}$. To simplify the model, we write $W_{1}$ as $W$, and we will try to determine the range of $\theta$ and $a^{*}, a_{*}$ such that the solution exists globally in time.

$\bullet$ The third goal of this paper is to consider weak linear noise effects associated with the phenomenon of wave breaking. Due to Theorem \ref{Global strong nonlinear} below, we see that if wave breaking occurs, the noise term does not grow
fast. Hence we consider $\theta=0$ in \eqref{SMCH21}, namely a non-autonomous pre-factor depending on time $t$. Precisely, we consider the MCH2 equation with linear multiplicative noise.
We will study the conditions that lead to the global existence and the blow-up in finite time of the solution, and then analyze the associated probabilities.

\section{Notation and preliminaries}

In this section, we begin by introducing some notations and recall some elementary results, for completeness, we list the lemmas and skip their some proofs for conciseness.

Let $L^{2}$ be the usual space of square-integrable functions on $\mathbb R$. For any real number $s\in \mathbb R$, $D^{s}= (1-\partial_{x}^{2})^{s/2}$ is defined by $\widehat{D^{s}f}(x)=(1+x^{2})^{\frac{s}{2}}\widehat{f}(x)$, where $\widehat{f}$ is the Fourier transform of $f$. The Sobolev space $H^{s}$ is defined as

\begin{equation*}\label{Hs}
H^{s}\triangleq \{f\in L^{2}(\mathbb R): \|f\|^{2}_{H^{s}}=\int_{\mathbb R}(1+x^{2})^{s}|\widehat{f}(x)|^{2}dx<\infty\},
\end{equation*}
and the inner product
\begin{equation*}\label{fg}
(f,g)_{H^{s}}:=\int_{\mathbb R}(1+x^{2})^{s}\widehat{f}(x)\overline{\widehat{g}}(x)dx=(\widehat{D^{s}f},\widehat{D^{s}g})_{L^{2}}.
\end{equation*}

In addition, $x\lesssim y$, $x,y\in\R$ means that there exists $C>0$, which may vary from line to line and depend on various parameters, such that $x\leq Cy$. Hereafter, $C$ denotes a positive constant, whose value may change from one place to another.

Firstly, we summarize some auxiliary results, which will be used to prove our main results. Define the regularizing operator $T_{\epsilon}$ on $\mathbb R$ as
\begin{eqnarray}\label{Tepsilon}
T_{\epsilon}f(x):=(1-\epsilon^{2}\Delta)^{-1}f(x)=\int_{\mathbb R}\frac{e^{i\xi x}\hat{f}(\xi)}{1+\epsilon^{2}|\xi|^{2}}d\xi, ~~\epsilon\in(0, 1).
\end{eqnarray}
Since $T_{\epsilon}$ can be characterized by its Fourier multipliers, see \cite{T}, it is easy to see that
\begin{eqnarray}\label{T}
&&[D^{s}, T_{\epsilon}]=0,\nonumber\\
&&(T_{\epsilon}f,g)_{L^{2}}=(f,T_{\epsilon}g)_{L^{2}},\nonumber\\
&&\|T_{\epsilon}u\|_{H^{s}}\leq \|u\|_{H^{s}}.
\end{eqnarray}
Where $[D^{s}, T_{\epsilon}]=D^{s}T_{\epsilon}-T_{\epsilon}D^{s}$. We therefore have the following lemma.
\begin{lem}\cite{T}\label{Te}
Let $f,g:\mathbb R\rightarrow\mathbb R$ such that $g\in W^{1,\infty}$ and $f\in L^{2}$. Then for some $C>0$,
\begin{eqnarray*}
\|[T_{\epsilon}, (g\cdot\nabla)]f\|_{L^{2}}\leq C\|g\|_{W^{1,\infty}}\|f\|_{L^{2}}.
\end{eqnarray*}
\end{lem}

Furthermore, we also need to  recall some useful commutator estimates.
\begin{lem}\cite{KP}\label{l1}
If $r>0$, then $H^{r}\bigcap L^{\infty}$ is an algebra. Moreover,
$\|uv\|_{H^{r}}\lesssim \|u\|_{L^{\infty}}\|v\|_{H^{r}}+\|u\|_{H^{r}}\|v\|_{L^{\infty}}$.
\end{lem}

\begin{lem}\cite{KP}\label{l2}
Let $r>0$, if $u\in H^{r}\bigcap W^{1,\infty}$ and $v\in H^{r-1}\bigcap L^{\infty}$, then
\begin{eqnarray*}
\|[D^{r},u]v\|_{L^{2}}\lesssim \|\partial_{x}u\|_{L^{\infty}}\|D^{r-1}v\|_{L^{2}}+\|D^{r}u\|_{L^{2}}\|v\|_{L^{\infty}},
\end{eqnarray*}
where $[D^{r},u]v=D^{r}uv-uD^{r}v$.
\end{lem}

A direct application of Lemma \ref{l1}-\ref{l2} gives the following estimates and we omit the
proof here.
\begin{lem}\label{F1,F2}
For the $F_{1}, F_{2}$ defined in \eqref{SMCH21} and for any $u,\gamma,u_{1},u_{2},\gamma_{1},\gamma_{2}\in H^{s}$ with $s>1/2$, we have
\begin{align*}
\|F_{1}(u,\gamma)\|_{H^{s}}
\lesssim& (\|u\|_{W^{1,\infty}}+\|\gamma\|_{W^{1,\infty}})(\|u\|_{H^{s}}+\|\gamma\|_{H^{s}}), \tag*{~$s>3/2,$}\nonumber\\
\|F_{2}(u,\gamma)\|_{H^{s}}
\lesssim& (\|u\|_{W^{1,\infty}}+\|\gamma\|_{W^{1,\infty}})(\|u\|_{H^{s}}+\|\gamma\|_{H^{s}}), \tag*{$s>3/2,$}\nonumber\\
\|F_{1}(u_{1},\gamma_{1})-F_{1}(u_{2},\gamma_{2})\|_{H^{s}}
\lesssim& (\|u_{1}\|_{H^{s}}+\|u_{2}\|_{H^{s}}+\|\gamma_{1}\|_{H^{s}}+\|\gamma_{2}\|_{H^{s}})\nonumber\\
&(\|u_{1}-u_{2}\|_{H^{s}}+\|\gamma_{1}-\gamma_{2}\|_{H^{s}}), \tag*{$s>3/2,$}\nonumber\\
\|F_{1}(u_{1},\gamma_{1})-F_{1}(u_{2},\gamma_{2})\|_{H^{s}}
\lesssim&(\|u_{1}\|_{H^{s+1}}+\|u_{2}\|_{H^{s+1}}+\|\gamma_{1}\|_{H^{s+1}}+\|\gamma_{2}\|_{H^{s+1}})\nonumber\\
&(\|u_{1}-u_{2}\|_{H^{s}}+\|\gamma_{1}-\gamma_{2}\|_{H^{s}}),  \tag*{$1/2<s<3/2,$}\nonumber\\
\|F_{2}(u_{1},\gamma_{1})-F_{2}(u_{2},\gamma_{2})\|_{H^{s}}
\lesssim&(\|u_{1}\|_{H^{s}}+\|\gamma_{2}\|_{H^{s}})(\|u_{1}-u_{2}\|_{H^{s}}+\|\gamma_{1}-\gamma_{2}\|_{H^{s}}),  \tag*{$s>3/2,$}\nonumber\\
\|F_{2}(u_{1},\gamma_{1})-F_{2}(u_{2},\gamma_{2})\|_{H^{s}}
\lesssim&(\|u_{1}\|_{H^{s+1}}+\|\gamma_{2}\|_{H^{s+1}})(\|u_{1}-u_{2}\|_{H^{s}}+\|\gamma_{1}-\gamma_{2}\|_{H^{s}}),  \tag*{$1/2<s<3/2.$}\nonumber\\
\end{align*}
\end{lem}
In addition, we provide the following algebraic inequality, which will be used in the proof of Theorem \ref{Global strong nonlinear}.
\begin{lem}\label{lem strong noise}
Let $c,M_{1}, M_{2}>0$. Assume $a,b_{*}, b^{*}>0$,
\begin{eqnarray*}
\text{either}~\eta>1,~0<\sqrt{b_{*}}<b(t)<\sqrt{b^{*}}~and~2b_{*}>b^{*}\\
\text{or}~ \eta=1, 0<\sqrt{b_{*}}<b(t)<\sqrt{b^{*}}~and ~2b_{*}>a+b^{*}.
\end{eqnarray*}
Then there is a constant $C>0$ such that for all $0\leq x_{1}\leq M_{1}y_{1}<\infty, 0\leq x_{2}\leq M_{2}y_{2}<\infty$,
{\small
\begin{align*}
\frac{a(x_{1}+x_{2})(y_{1}^{2}+y_{2}^{2})+b(t)(1+x_{1}+x_{2})^{\eta}(y_{1}^{2}+y_{2}^{2})}{1+y_{1}^{2}+y_{2}^{2}}
-&\frac{2b(t)(1+x_{1}+x_{2})^{\eta}(y_{1}^{2}+y_{2}^{2})^{2}}{(1+y_{1}^{2}+y_{2}^{2})^{2}}\nonumber\\
+&\frac{cb(t)(1+x_{1}+x_{2})^{\eta}(y_{1}^{2}+y_{2}^{2})^{2}}{(1+y_{1}^{2}+y_{2}^{2})^{2}(1+\log(1+y_{1}^{2}+y_{2}^{2}))}\leq C.
\end{align*}}
\end{lem}
\begin{pf}
Since $0\leq \frac{x_{1}}{M_{1}}\leq y_{1}<\infty, 0\leq \frac{x_{2}}{M_{2}}\leq y_{2}<\infty$, we obtain
{
\begin{align*}
&\frac{a(x_{1}+x_{2})(y_{1}^{2}+y_{2}^{2})+b(t)(1+x_{1}+x_{2})^{\eta}(y_{1}^{2}+y_{2}^{2})}{1+y_{1}^{2}+y_{2}^{2}}
-\frac{2b(t)(1+x_{1}+x_{2})^{\eta}(y_{1}^{2}+y_{2}^{2})^{2}}{(1+y_{1}^{2}+y_{2}^{2})^{2}}\nonumber\\
&+\frac{cb(t)(1+x_{1}+x_{2})^{\eta}(y_{1}^{2}+y_{2}^{2})^{2}}{(1+y_{1}^{2}+y_{2}^{2})^{2}(1+\log(1+y_{1}^{2}+y_{2}^{2}))}\nonumber\\
\leq& a(x_{1}+x_{2})+b^{*}(1+x_{1}+x_{2})^{\eta}-2b_{*}(1+x_{1}+x_{2})^{\eta}\frac{(y_{1}^{2}+y_{2}^{2})^{2}}{(1+y_{1}^{2}+y_{2}^{2})^{2}}
+\frac{cb^{*}(1+x_{1}+x_{2})^{\eta}}{1+\log(1+(\frac{x_{1}}{M_{1}})^{2}+(\frac{x_{2}}{M_{2}})^{2})}.
\end{align*}}
When $\eta>1$ and $2b_{*}>b^{*}$ or $\eta=1$ and $2b_{*}>a+b^{*}$, we find that the inequality will tend to $-\infty$ for $x_{1}\rightarrow\infty$ or $x_{2}\rightarrow\infty$ (namely $y_{1}\rightarrow\infty$ or $y_{2}\rightarrow\infty$), which completes the proof.

\end{pf}

Finally, we present the following lemma to establish the Theorem \ref{noise I} on the global existence solutions.
\begin{lem}\label{process control}
Let $\alpha(t)$ be a deterministic and locally bounded function. Suppose that $\lambda>0$,  and $x(t)$ satisfies
\begin{eqnarray*}
x(t)=e^{\int_{0}^{t}\alpha(t^{'})d W_{t^{'}}-\int_{0}^{t}\lambda\alpha^2(t^{'})d t^{'}}.
\end{eqnarray*}
For $R>1$ define $\tau_{R}=\inf\{t\geq0: x(t)>R\}$. Then we have
\begin{eqnarray*}
\mathbb P\{\tau_{R}=\infty\}\geq1-R^{-2\lambda}.
\end{eqnarray*}
\end{lem}
\begin{pf}
Noting that $$x(t)^{2\lambda}=e^{\int_{0}^{t}2\lambda\alpha(t^{'})d W_{t^{'}}-\frac{1}{2}\int_{0}^{t}(2\lambda)^2\alpha^2(t^{'})d t^{'}}$$ is an exponential martingale, by the martingale stopping theorem, we can derive
$\mathbb{E}x(t\wedge\tau_{R})^{2\lambda}=1$. This yields
\begin{eqnarray*}
\mathbb P(\tau_{R}=\infty)=\lim_{n\rightarrow\infty}\mathbb P(\tau_{R}>n)=\lim_{n\rightarrow\infty}\mathbb P(x(n\wedge\tau_{R})^{2\lambda}<R^{2\lambda})\geq\lim_{n\rightarrow\infty}\bigg(1-\frac{\mathbb Ex(n\wedge\tau_{R})^{2\lambda}}{R^{2\lambda}}\bigg)=1-R^{-2\lambda},
\end{eqnarray*}
which finishes the proof.
\end{pf}
\begin{lem}\label{energy estimation}
Let $s>3/2$, $F_{1}, F_{2}$ and $T_{\epsilon}$ be given in \eqref{SMCH21} and \eqref{Tepsilon} respectively. Then there is a constant $K=K(s)>0$ such that for all $\epsilon>0$,
\begin{eqnarray*}
&&|(T_{\epsilon}[uu_{x}],T_{\epsilon}u)_{H^{s}}|+|(T_{\epsilon}F_{1}(u,\gamma),T_{\epsilon}u)_{H^{s}}|+|(T_{\epsilon}[u\gamma_{x}],T_{\epsilon}\gamma)_{H^{s}}|+|(T_{\epsilon}F_{2}(u,\gamma),T_{\epsilon}\gamma)_{H^{s}}|\nonumber\\
&\leq& K(\|u\|_{W^{1,\infty}}+\|\gamma\|_{W^{1,\infty}})(\|u\|_{H^{s}}^{2}+\|\gamma\|_{H^{s}}^{2}).
\end{eqnarray*}
\begin{pf}
According to \eqref{T}, we derive
\begin{eqnarray*}
(T_{\epsilon}[uu_{x}],T_{\epsilon}u)_{H^{s}}&=&(D^{s}T_{\epsilon}[uu_{x}],D^{s}T_{\epsilon}u)_{L^{2}}\nonumber\\
&=&([D^{s},u]u_{x}, D^{s}T_{\epsilon}^{2}u)_{L^{2}}+([T_{\epsilon},u]D^{s}u_{x}, D^{s}T_{\epsilon}u)_{L^{2}}+(uD^{s}T_{\epsilon}u_{x}, D^{s}T_{\epsilon}u)_{L^{2}},\nonumber\\
(T_{\epsilon}[u\gamma_{x}],T_{\epsilon}\gamma)_{H^{s}}&=&(D^{s}T_{\epsilon}[u\gamma_{x}],D^{s}T_{\epsilon}\gamma)_{L^{2}}\nonumber\\
&=&([D^{s},u]\gamma_{x}, D^{s}T_{\epsilon}^{2}\gamma)_{L^{2}}+([T_{\epsilon},u]D^{s}\gamma_{x}, D^{s}T_{\epsilon}\gamma)_{L^{2}}+(uD^{s}T_{\epsilon}\gamma_{x}, D^{s}T_{\epsilon}\gamma)_{L^{2}}.
\end{eqnarray*}
Then by Lemma \ref{Te}, Lemma \ref{l1}, \eqref{T} and Sobolev embedding Theorem, we have
\begin{eqnarray*}
|(T_{\epsilon}[uu_{x}],T_{\epsilon}u)_{H^{s}}|+|(T_{\epsilon}[u\gamma_{x}],T_{\epsilon}\gamma)_{H^{s}}|\lesssim (\|u\|_{W^{1,\infty}}+\|\gamma\|_{W^{1,\infty}})(\|u\|_{H^{s}}^{2}+\|\gamma\|_{H^{s}}^{2}).
\end{eqnarray*}
In addition, using Lemma \ref{F1,F2} and \eqref{T}, we obtain that
\begin{eqnarray*}
|(T_{\epsilon}F_{1}(u,\gamma),T_{\epsilon}u)_{H^{s}}|+|(T_{\epsilon}F_{2}(u,\gamma),T_{\epsilon}\gamma)_{H^{s}}|\lesssim (\|u\|_{W^{1,\infty}}+\|\gamma\|_{W^{1,\infty}})(\|u\|_{H^{s}}^{2}+\|\gamma\|_{H^{s}}^{2}).
\end{eqnarray*}
Combining the above two inequalities, we complete the proof.
\end{pf}
\end{lem}

\section{Local existence and uniqueness for SMCH2}

In this section, we consider the following stochastic system with multiplicative noise \eqref{SMCH21} in $H^{s}(\R)\times H^{s}(\R)$.

\subsection{Assumptions}
For the main results in this paper, we rely on the following different assumptions concerning random perturbation term in \eqref{SMCH21}.
We assume that $(h_{1}, h_{2}): [0, \infty)\times (H^{s}\times H^{s})\ni (t,u,\gamma)\rightarrow (h_{1}(t,u,\gamma), h_{2}(t,u,\gamma))\in L_{2}(U, H^{s}\times H^{s})$ are continuous in $(t,u,\gamma)$. Moreover, we assume
\begin{asu}\label{h1,h2}
(H.1) There exists some non-decreasing function $f:[0,\infty)\to[0,\infty)$ with $f(0)=0$ such that for all $(u, \gamma)\in H^{s}\times H^{s}, s>1/2$,
\begin{eqnarray}\label{h1}
\sum_{i=1}^2\|h_{i}(t, u,\gamma)\|_{L_{2}(U, H^{s})}\leq f(\|u\|_{W^{1,\infty}}+\|\gamma\|_{W^{1,\infty}})(1+\|u\|_{H^{s}}+\|\gamma\|_{H^{s}}).
\end{eqnarray}
(H.2) There exists some non-decreasing function $g:[0,\infty)\to[0,\infty)$ such that for all $(u_{1},\gamma_{1}), (u_{2},\gamma_{2})\in H^{s}\times H^{s}, s>1/2$,
\begin{eqnarray*}
\sup_{\|u_{1}\|_{H^{s}},\|\gamma_{1}\|_{H^{s}},\|u_{2}\|_{H^{s}},\|\gamma_{2}\|_{H^{s}}\leq N}&&\sum_{i=1}^2\|h_{i}(t, u_{1},\gamma_{1})-h_{i}(t, u_{2},\gamma_{2})\|_{L_{2}(U, H^{s})}\nonumber\\
&\leq& g(N)\cdot(\|u_{1}-u_{2}\|_{H^{s}}+\|\gamma_{1}-\gamma_{2}\|_{H^{s}}), N\geq1.
\end{eqnarray*}
\end{asu}
\begin{asu}\label{nonLinear noise} $h_{1}(t, u,\gamma)d\mathcal{W}_{1}=a(t)(1 +\|u\|_{W^{1,\infty}}+\|\gamma\|_{W^{1,\infty}})^{\theta}udW, h_{2}(t, u,\gamma)d\mathcal{W}_{2}=a(t)(1+\|u\|_{W^{1,\infty}}+\|\gamma\|_{W^{1,\infty}})^{\theta}\gamma dW$ for a standard 1-D Brownian motion $W$ and $\theta>0$, $0<a_{\ast}\leq a^{2}(t)\leq a^{\ast}$ for all $t$.
\end{asu}

\begin{asu}\label{Linear noise}
$h_{1}(u,\gamma)d\mathcal{W}_{1}=b(t)udW$, $h_{2}(u,\gamma)d\mathcal{W}_{2}=b(t)\gamma dW$ for a standard 1-D Brownian motion $W$, and there are constants $b_{\ast}, b^{\ast}>0$ such that $0<b_{\ast}\leq b^{2}(t)\leq b^{\ast}$ for all $t$.
\end{asu}
%
\subsection{Definitions of the solutions.}Next, we give the definition of pathwise solution to \eqref{SMCH21}.
\begin{defn}\label{Pathwise}
(Pathwise solutions). Let $\mathcal{S}=(\Omega,\mathcal{F},\mathbb{P}, \{\mathcal{F}_{t}\}_{t\geq0},\mathcal{W}_{1}, \mathcal{W}_{2})$ be a fixed stochastic basis. Let $s>3/2$ and $z_{0}=(u_{0},\gamma_{0})$ be an $H^{s}\times H^{s}$-valued $\mathcal{F}_{0}$-measurable random variable.\\
\text{1}.A local pathwise solution to \eqref{SMCH21} is a pair $(z,\tau)$, where $\tau\geq0$ is a stopping time satisfying $\mathbb P\{\tau>0\}=1$ and $z=(u, \gamma):\Omega\times[0,\tau)\rightarrow H^{s}\times H^{s}$ is an $\mathcal{F}_{t}$-adapted $H^{s}\times H^{s}$-valued process satisfying $\mathbb{P}-a.s.$
\begin{eqnarray}\label{z}
z\in C([0,\tau); H^{s}\times H^{s}),
\end{eqnarray}
and $\mathbb{P}-a.s.$,
\begin{eqnarray*}\label{z1}
&&u(t)-u(0)+\int_{0}^{t}[uu_{x}+F_{1}(u,\gamma)]d t^{'}=\int_{0}^{t}h_{1}(t^{'}, u,\gamma)d\mathcal{W}_{1},\\
&&\gamma(t)-\gamma(0)+\int_{0}^{t}[u\gamma_{x}+F_{2}(u,\gamma)]d t^{'}=\int_{0}^{t}h_{2}(t^{'},u,\gamma)d\mathcal{W}_{2},\ \ t\in [0,\tau).
\end{eqnarray*}
\text{2}. Local pathwise uniqueness: if given any two local pathwise solutions $(z_{1}, \tau_{1})$ and $(z_{2}, \tau_{2})$ with $\mathbb P\{z_{1}(0)=z_{2}(0)\}=1$, we have
\begin{eqnarray*}\label{Z1=Z2}
\mathbb P\{z_{1}(t)=z_{2}(t), \ \ t\in [0, \tau_{1}\wedge \tau_{2})\}=1.
\end{eqnarray*}
\text{3}.Additionally, $(z, \tau^{\ast})$ is called a maximal pathwise solution to \eqref{SMCH21} if $\tau^{\ast}>0$ almost
surely and there is an increasing sequence $\tau_{n}\rightarrow\tau^{\ast}$ such that for any $n\in \mathbb N$, $(z, \tau_{n})$ is a pathwise solution to \eqref{SMCH21} and on the set $\{\tau^{\ast}<\infty\}$,
\begin{eqnarray*}\label{t}
\mathop{\sup}\limits_{t\in[0,\tau_{n}]}(\|u\|_{H^{s}}+\|\gamma\|_{H^{s}})\geq n, ~n\geq1.
\end{eqnarray*}
\text{4}. If  $(z, \tau^{\ast})$ is a maximal pathwise solution and $\tau^{\ast}=\infty$ almost surely, then we call that the pathwise solution exists globally.
\end{defn}

\subsection{Main results and remarks.}
Now, we summarize our major contributions, such as existence of pathwise solutions, global well-posedness of \eqref{SMCH21} and the blow-up results, and the concrete proofs will be provided later in the remainder of the paper.
\begin{thm}\label{Maximal solutions}(Maximal solutions)
Let $s>3/2$, and $h_{1}(t, u,\gamma), h_{2}(t, u,\gamma)$ satisfy Assumption 3.1. For a given stochastic basis $\mathcal{S}=(\Omega, \mathcal{F}, \mathbb P, \{\mathcal{F}_{t}\}_{t\geq0},\mathcal{W}_{1}, \mathcal{W}_{2})$, if $(u_{0},\gamma_{0})$ is an $H^{s}\times H^{s}$-valued $\mathcal{F}_{0}$-measurable random variable, then there is a local unique pathwise solution $(z,\tau)$ to \eqref{SMCH21} in
the sense of Definition \ref{Pathwise} with
\begin{eqnarray*}\label{Local existence}
z\in C([0,\tau); H^{s}\times H^{s}).
\end{eqnarray*}
Moreover, $(z,\tau)$ can be extended to a unique maximal pathwise solution $(z,\tau^{\ast})$ and the following blow up scenario satisfies $\mathbb P-a.s.$ on the set $\{\tau^\ast<\infty\}$,
\begin{eqnarray}\label{Blow-up criterion}
&&1_{\{\lim\mathop{\sup}\limits_{t\rightarrow \tau^{\ast}}(\|u(t)\|_{H^{s}}+\|\gamma(t)\|_{H^{s}})=\infty\}}=1_{\{\lim\mathop{\sup}\limits_{t\rightarrow \tau^{\ast}}(\|u(t)\|_{W^{1,\infty}}+\|\gamma(t)\|_{W^{1,\infty}})=\infty\}}.
\end{eqnarray}
\end{thm}
\begin{re}
The proof of Theorem \ref{Maximal solutions} combines the techniques as used in the papers \cite{T,NV,ANR,DEM,DEM1,DM}. By constructing the approximate sequence of the truncation problem of $W^{1,\infty}\times W^{1,\infty}$. Such a cut-off means linear growth of $u$ and $\gamma$, and guarantees the global existence of an approximate solution.

\end{re}

Turning to noise-driven regularization effects, the blow-up scenario \eqref{Blow-up criterion} suggests relating the noise coefficient to the $W^{1,\infty}\times W^{1,\infty}$ of $(u, \gamma)$. Therefore we consider scalable noise impact, i.e. we assume $h_{1}(t,u,\gamma)d\mathcal{W}_{1}=a(t)(1+\|u\|_{W^{1,\infty}}+\|\gamma\|_{W^{1,\infty}})^{\theta}u dW, h_{2}(t,u,\gamma)d\mathcal{W}_{2}=a(t)(1+\|u\|_{W^{1,\infty}}+\|\gamma\|_{W^{1,\infty}})^{\theta}\gamma dW$ for a standard 1-D Brownian motion $W$, some $\theta>0$, $0<a_{\ast}\leq a^{2}(t)\leq a^{\ast}$ for all $t$. When $a^{*}, a_{*}$ and $\theta$ satisfy certain stronger conditions, the noise term remove the formation of singularities.
\begin{thm}\label{Global strong nonlinear}
(Global existence for strong nonlinear noise). Let Assumption \ref{nonLinear noise} hold and assume that $\mathcal{S}=(\Omega,\mathcal{F},\mathbb{P},\{\mathcal{F}_{t}\}_{t\geq0},W)$ is a fixed stochastic basis. Let $s>5/2$, $(u_{0}, \gamma_{0})\in H^{s}\times H^{s}$ be an $H^{s}\times H^{s}$-valued $\mathcal{F}_{0}$-measurable random variable. Assume that $\theta$ and $a^{*}, a_{*}(0<a_{\ast}\leq a^{2}(t)\leq a^{\ast})$ satisfy
\begin{eqnarray*}\label{thm}
either~2a_{*}>a^{*},~\theta>1/2~or~2a_{*}>K+a^{*}, ~\theta=1/2,
\end{eqnarray*}
where $K=K(s)$ is the constant introduced in Lemma \ref{energy estimation}. Then the corresponding maximal solution $(z, \tau^{\ast})$ to \eqref{SMCH21} satisfies
\begin{eqnarray*}
\mathbb P\{\tau^{\ast}=\infty\}=1.
\end{eqnarray*}
\end{thm}
\begin{re}
Theorem \ref{Global strong nonlinear} means that blow-up of pathwise solutions might only be observed if the noise is weak. According to Theorem \ref{Global strong nonlinear}, we can see that if wave breaking occurs, the noise term will not bring rapid growth. Therefore, we consider $\theta=0$ but  a non-autonomous pre-factor dependent on time $t$ is introduced.
\end{re}

To detect such noise, we analyze the simpler form $h_{1}(t,u,\gamma)d\mathcal{W}_{1}= b(t)u d W, h_{2}(t,u,\gamma)d\mathcal{W}_{2} = b(t)\gamma d W$, $W$ is a standard 1-D Brownian motion. Even in this linear noise case the situation is quite interesting allowing for global existence as well as blow-up of solutions. For global existence, we can identify two cases.

Using Lemma \ref{l1}-Lemma \ref{F1,F2} and the integration by parts, we conclude that there is a $C=C(s)>1$ such that
\begin{align}\label{PPY}
&-\int_{\mathbb R}D^{s}v_{1}D^{s}(v_{1}v_{1x})dx-\int_{\mathbb R}D^{s}v_{1}D^{s}F_{1}(v_{1},v_{2})dx-\int_{\mathbb R}D^{s}v_{2}D^{s}(v_{1}v_{2x})dx-\int_{\mathbb R}D^{s}v_{2}D^{s}F_{2}(v_{1},v_{2})dx\nonumber\\
&\leq \frac{1}{2} C(\|v_{1}\|_{W^{1,\infty}}+\|v_{2}\|_{W^{1,\infty}})(\|v_{1}\|_{H^{s}}^{2}+\|v_{2}\|_{H^{s}}^{2}).
\end{align}

\begin{thm}\label{noise I}
(Global existence for weak noise I). Let $s>3/2$, Assumption \ref{Linear noise} be verified and $\mathcal{S}=(\Omega,\mathcal{F},\mathbb{P},\{\mathcal{F}_{t}\}_{t\geq0},W)$ be a fixed stochastic basis. Assume $(u_{0}, \gamma_{0})$ is an $H^{s}\times H^{s}$-valued $\mathcal{F}_{0}$-measurable random variable. Let $Q=Q(s)>0$ be the constant such that
the embedding $\|u\|_{W^{1,\infty}}<Q\|u\|_{H^{s}}, \|\gamma\|_{W^{1,\infty}}<Q\|\gamma\|_{H^{s}}$ holds. Let $C=C(s)>1$ be in \eqref{PPY}. If there is a  $R>1$ and
$\lambda_{1}>1$ satisfying $\mathbb P-a.s.$
\begin{eqnarray*}
\|u_{0}\|_{H^{s}}^{2}+\|\gamma_{0}\|_{H^{s}}^{2}<\frac{b_{\ast}^{2}}{4C^{2}Q^{2}\lambda_{1}^{2}R},
\end{eqnarray*}
then \eqref{SMCH21} has a maximal solution $(z,\tau^{\ast})$ satisfying for any $0<\lambda_{2}<\frac{\lambda_{1}-1}{\lambda_{1}}$ the estimate
\begin{eqnarray*}
\mathbb P\bigg\{ \|u(t)\|_{H^{s}}^{2}+\|\gamma(t)\|_{H^{s}}^{2}<\frac{b_{\ast}^{2}}{C^{2}Q^{2}\lambda_{1}^{2}}~~for ~all~ t>0\bigg\}\geq 1-R^{-2\lambda_{2}}.
\end{eqnarray*}
\begin{re}
Theorem \ref{noise I} presents a global existence solution with bounded initial data. This result can not be observed in the deterministic case.
\end{re}

\end{thm}

\begin{thm}\label{noise II}
(Global existence for weak noise II) Let $s>5/2$, Assumption \ref{Linear noise} be verified and $\mathcal{S}=(\Omega,\mathcal{F},\mathbb{P},\{\mathcal{F}_{t}\}_{t\geq0},W)$ be a fixed stochastic basis. $(u_{0}, \gamma_{0})$ is an $H^{s}\times H^{s}$-valued $\mathcal{F}_{0}$-measurable random variable. If
\begin{eqnarray*}
&&\mathbb P\{(1-\partial_{x}^{2})u_{0}(x)>0, \forall x\in\mathbb R\}=p, ~~\mathbb P\{(1-\partial_{x}^{2})u_{0}(x)<0, \forall x\in\mathbb R\}=q,
\end{eqnarray*}
and  there exists some $x_{0}\in\mathbb R$ such that
\begin{eqnarray*}
&&\mathbb P\{(1-\partial_{x}^{2})u_{0}(x)\leq0, ~x\leq x_{0} ~~and~~ (1-\partial_{x}^{2})u_{0}(x)\geq0, ~x\geq x_{0}\}=m,
\end{eqnarray*}

for some $p, q, m\in[0, 1]$, then the corresponding maximal solution $(z, \tau^{\ast})$ to \eqref{SMCH21} satisfies
\begin{eqnarray*}
\mathbb P\{\tau^{\ast}=\infty\}\geq p+q+m.
\end{eqnarray*}
\end{thm}
\begin{re}
The proof of Theorem \ref{noise II} depends on the analysis of a PDE with random coefficient. When $b(t)=0$ and taking $(p, q, m)=(1, 0, 0)$, $(p, q, m)=(0, 1, 0)$ or $(p, q, m)=(0, 0, 1)$ in Theorem \ref{noise II}, we obtain the global existence for the deterministic MCH2 system. Therefore, in this sense,
Theorem \ref{noise II} covers the deterministic result.
\end{re}

\begin{thm}\label{blow up initial value}(Wave breaking criterion for weak noise I)
Let $\mathcal{S}=(\Omega,\mathcal{F},\mathbb{P},\{\mathcal{F}_{t}\}_{t\geq0}, W)$ be a fixed stochastic basis and $s>5/2$. Let Assumption \ref{Linear noise} be verified and $(u_{0}, \gamma_{0})$ be an $H^{s}\times H^{s}$-valued $\mathcal{F}_{0}$-measurable random variable. If for some $c\in(0,1)$ and $x_{0}\in\mathbb R$,
\begin{eqnarray}\label{u0x}
u_{0x}(x_{0})<-\frac{1}{2}\sqrt{\frac{(b^{*})^{2}}{c^{2}}+4(\|u_{0}\|_{H^{1}}^{2}+\|\gamma_{0}\|_{H^{1}}^{2})}-\frac{b^{*}}{2c}~\mathbb P-a.s.,
\end{eqnarray}
then the maximal solution $(z, \tau^{\ast})$ to \eqref{SMCH21} satisfies
\begin{eqnarray*}
\mathbb P\{\tau^{*}<\infty\}\geq\mathbb P\left\{e^{\int_{0}^{t}b(t^{'})dW_{t^{'}}+\int_{0}^{t}\frac{b^{*}-b^{2}(t^{'})}{2}dt^{'}}\geq c ~for~ all~ t\right\}>0.
\end{eqnarray*}
\end{thm}
\begin{re}
Theorem \ref{blow up initial value} detects the solution singularities in finite time under certain initial data, while Theorem \ref{noise I} provides a global existence result. We stress that these two results do not contain each other. In Theorem \ref{noise I}, assuming $\|u_{0}\|_{H^{s}}^{2}+\|\gamma_{0}\|_{H^{s}}^{2}<\frac{b_{\ast}^{2}}{2C^{2}Q^{2}\lambda_{1}^{2}R}$, then $z$ globally exists with probability greater than $1-R^{-2\lambda_{2}}$. In Theorem \ref{blow up initial value}, \eqref{u0x} implies that $\|u_{0}\|_{H^{s}}^{2}>\frac{1}{Q^{2}}\|u_{0}\|_{W^{1,\infty}}^{2}>\frac{(b^{\ast})^{2}}{c^{2}Q^{2}}>\frac{b_{\ast}^{2}}{2C^{2}Q^{2}\lambda_{1}^{2}R}$.

\end{re}

\begin{thm}\label{blow up initial value 2}(Wave breaking criterion for weak noise II)
Let $\mathcal{S}=(\Omega,\mathcal{F},\mathbb{P},\{\mathcal{F}_{t}\}_{t\geq0}, W)$ be a fixed stochastic basis and $s>5/2$. Let Assumption \ref{Linear noise} be verified and $(u_{0}, \gamma_{0})$ be an $H^{s}\times H^{s}$-valued $\mathcal{F}_{0}$-measurable random variable. If for some $c\in(0,1)$,
\begin{align}\label{u0x1}
\int_{\mathbb R}u_{0x}^{3}(x)dx&<-\sqrt{\frac{(b^{*})^{2}}{4c^{2}}(\|u_{0}\|_{H^{1}}^{2}+\|\gamma_{0}\|_{H^{1}}^{2})^{2}+\frac{15}{8} (\|u_{0}\|_{H^{1}}^{2}+\|\gamma_{0}\|_{H^{1}}^{2})^{3}}\nonumber\\
&\qquad\qquad\quad-\frac{b^{*}}{2c}(\|u_{0}\|_{H^{1}}^{2}+\|\gamma_{0}\|_{H^{1}}^{2})~\mathbb P-a.s.,\nonumber\\
\end{align}
then the maximal solution $(z, \tau^{\ast})$ to \eqref{SMCH21} satisfies
\begin{eqnarray*}
\mathbb P\{\tau^{*}<\infty\}\geq\mathbb P\left\{e^{\int_{0}^{t}b(t^{'})dW_{t^{'}}+\int_{0}^{t}\frac{b^{*}-b^{2}(t^{'})}{2}dt^{'}}\geq c ~for~ all~ t\right\}>0.
\end{eqnarray*}
\end{thm}

\begin{re}
Theorem \ref{blow up initial value 2} detects the solution singularities in finite time under certain initial data. \eqref{u0x1} implies that $\int_{\mathbb R}u_{0x}^{3}(x)dx<-\frac{b^{*}}{c}(\|u_{0}\|_{H^{1}}^{2}+\|\gamma_{0}\|_{H^{1}}^{2})$, which combined with $min_{x\in\mathbb R}u_{0x}(x)(\|u_{0}\|_{H^{1}}^{2}+\|\gamma_{0}\|_{H^{1}}^{2})\leq\int_{\mathbb R}u_{0x}^{3}(x)dx$ derives $\|u_{0}\|_{H^{s}}^{2}>\frac{1}{Q^{2}}\|u_{0}\|_{W^{1,\infty}}^{2}>\frac{(b^{\ast})^{2}}{c^{2}Q^{2}}>\frac{b_{\ast}^{2}}{2C^{2}Q^{2}\lambda_{1}^{2}R}$. So, the initial value conditions here and those given in Theorem \ref{noise I} do not contain each other.

\end{re}

\section{Sketch of the Proof of Theorem \ref{Maximal solutions}}
We consider the initial value problem \eqref{SMCH21}. The proof of existence and uniqueness of pathwise solutions can be carried out by standard procedures used in many works, see \cite{T,NV,DEM,DEM1,GZ,CH} for more details. Therefore we only give a sketch.\\
$1.$ (Approximation scheme) The first step is to construct a suitable approximation scheme. For any $R>1$, we let $\chi_{R}(x):[0,\infty)\rightarrow[0,1]$ be a $C_0^{\infty}$ function such that $\chi_{R}(x)=1$ for $x\in[0,R]$ and $\chi_{R}(x)=0$ for $x>2R$. Then we consider the following cut-off problem on $\mathbb{R}$,
\begin{equation}\label{SMCH2 cut-off}
\left\{
\begin{array} {l}
d u+\chi_{R}(\|u\|_{W^{1,\infty}}+\|\gamma\|_{W^{1,\infty}})[uu_x+F_{1}(u,\gamma)]d t=\chi_{R}(\|u\|_{W^{1,\infty}}+\|\gamma\|_{W^{1,\infty}})h_{1}(t, u,\gamma)d\mathcal{W}_{1},~t>0, \\
d \gamma+\chi_{R}(\|u\|_{W^{1,\infty}}+\|\gamma\|_{W^{1,\infty}})[u\gamma_x+F_{2}(u,\gamma)]d t=\chi_{R}(\|u\|_{W^{1,\infty}}+\|\gamma\|_{W^{1,\infty}})h_{2}(t, u,\gamma)d\mathcal{W}_{2},~t>0, \\
u(\omega,0,x)=u_0(\omega,x),\\
\gamma(\omega,0,x)=\gamma_0(\omega,x).
\end{array}
\right.
\end{equation}
From \eqref{F1,F2}, we observe that the nonlinear term $F_{1}(u,\gamma), F_{2}(u,\gamma),$ preserves the $H^{s}\times H^{s}$-regularity of $(u,\gamma)$ for any $s>3/2$. However, in order to apply the stochastic differential equation (SDE) theory in Hilbert space to \eqref{SMCH2 cut-off}, we will mollify the transport term $uu_x, u\gamma_x$ since the products $uu_x$ and $u\gamma_x$ lose one regularity. For this reason, we consider the following approximation scheme:
\begin{equation}\label{SMCH2 mollify}
\left\{
\begin{array} {l}
d u+G_{1,\epsilon}(u,\gamma)d t=\chi_{R}(\|u\|_{W^{1,\infty}}+\|\gamma\|_{W^{1,\infty}})h_{1}(t, u,\gamma)d\mathcal{W}_{1},~t>0,~x\in\R, \\
d \gamma+G_{2,\epsilon}(u,\gamma)d t=\chi_{R}(\|u\|_{W^{1,\infty}}+\|\gamma\|_{W^{1,\infty}})h_{2}(t, u,\gamma)d\mathcal{W}_{2},~t>0,~x\in\R, \\
G_{1,\epsilon}(u,\gamma)=\chi_{R}(\|u\|_{W^{1,\infty}}+\|\gamma\|_{W^{1,\infty}})[J_{\epsilon}((J_{\epsilon}u)(J_{\epsilon}u)_x)+F_{1}(u,\gamma)],\\
G_{2,\epsilon}(u,\gamma)=\chi_{R}(\|u\|_{W^{1,\infty}}+\|\gamma\|_{W^{1,\infty}})[J_{\epsilon}((J_{\epsilon}u)(J_{\epsilon}\gamma)_x)+F_{2}(u,\gamma)],\\
u(0,x)=u_0(x)\in H^{s},~\gamma(0,x)=\gamma_0(x)\in H^{s},
\end{array}
\right.
\end{equation}
where $J_{\epsilon}$ is the Friedrichs mollifier. According to the theory of SDE in Hilbert space (see for
example \cite{PR,KX}), for a fixed stochastic basis $\mathcal{S}=(\Omega, \mathcal{F}, \mathbb P, \{\mathcal{F}_{t}\}_{t\geq0},\mathcal{W}_{1}, \mathcal{W}_{2})$ and for $(u_0,\gamma_0)\in  H^{s}\times H^{s}$ with $s>5/2$,  \eqref{SMCH2 mollify} admits a unique solution $(u_{\epsilon}, \gamma_{\epsilon})\in C([0, T_{\epsilon}), H^{s}\times H^{s})$.

In addition, the uniform $L^{\infty}(\Omega;W^{1,\infty}\times W^{1,\infty})$ condition provided by the cut-off function $\chi_{R}$ enables us to split the expectation $\mathbb E(\|u_{\epsilon}\|^{2}_{H^{s}}\|u_{\epsilon}\|_{W^{1,\infty}}|\mathcal{F}_0)$, $\mathbb E(\|u_{\epsilon}\|^{2}_{H^{s}}\|\gamma_{\epsilon}\|_{W^{1,\infty}}|\mathcal{F}_0)$ to close a priori $L^{2}(\Omega, H^{s}\times H^{s}, \mathbb{P}(\cdot|\mathcal{F}_0))$ estimate for $u_{\epsilon}, \gamma_{\epsilon}$. Then we can go along the lines as we prove Lemma \ref{blp} to find that for each fixed $\epsilon$, if $T_{\epsilon}<\infty$, then $\lim\sup_{t\rightarrow T_{\epsilon}}(\|u_{\epsilon}\|_{W^{1,\infty}}+\|\gamma_{\epsilon}\|_{W^{1,\infty}})=\infty$. Due to the cut-off in \eqref{SMCH2 mollify} for $\mathbb{P}$-a.s. $\omega\in\Omega$, $\|u_{\epsilon}\|_{W^{1,\infty}}, \|\gamma_{\epsilon}\|_{W^{1,\infty}}$ are always bounded and hence $(u_{\epsilon}, \gamma_{\epsilon})$ is actually a global in time solution, that is $(u_{\epsilon}, \gamma_{\epsilon})\in C([0,\infty),H^{s}\times H^{s})~ \mathbb P-a.s.$\\
$2.$ (Pathwise solution to the cut-off problem in $H^{s}\times H^{s}$ with $s>5/2$) By applying the stochastic compactness arguments from Prokhorov's and Skorokhod's Theorem, we obtain the almost sure convergence for a new approximation solution $((\tilde{u_{\epsilon}}, \tilde{\gamma_{\epsilon}}), \tilde{\mathcal{W}_{1\epsilon}}, \tilde{\mathcal{W}_{2\epsilon}})$ defined on a new probability space. By virtue of a refined martingale representation Theorem \cite[Theorem A.1]{H1}, we may set $\epsilon\rightarrow0$ in $((\tilde{u_{\epsilon}}, \tilde{\gamma_{\epsilon}}), \tilde{\mathcal{W}_{1\epsilon}}, \tilde{\mathcal{W}_{2\epsilon}})$ to obtain a martingale solution in $H^{s}\times H^{s}$ with $s>5/2$. Here, the Gy\"{o}ngy-Krylov characterization \cite{GK} of the convergence in probability can be used here to prove the convergence of the original approximation solutions, and one can refer to \cite[Theorem 1.7]{T} for more details.

Finally, since $F_{1}(u,\gamma), F_{2}(u,\gamma)$ satisfy the estimates as in Lemma \ref{F1,F2} and $h_{1}(t, u,\gamma), h_{2}(t, u,\gamma)$ satisfies Assumption \ref{h1,h2}, we conclude that $G_{1,\epsilon}, G_{2,\epsilon}, \chi_{R}(\|u\|_{W^{1,\infty}}+\|\gamma\|_{W^{1,\infty}})h_1, \chi_{R}(\|u\|_{W^{1,\infty}}+\|\gamma\|_{W^{1,\infty}})h_2$ are Lipschitz continuous. So, one can obtain the pathwise uniqueness easily. Then by the Yamada-Watanabe principle, we derive the existence and uniqueness of the pathwise solution to \eqref{SMCH2 mollify} denoted by $z^{R}=(u^R,\gamma^R)$. \\
$3.$ (Remove the cut-off and extend the range of $s$ to $s>3/2$)

Let $\tau_R:=R\wedge\inf\{t\geq 0: \|u^R(t)\|_{W^{1,\infty}}+\|\gamma^R(t)\|_{W^{1,\infty}}>R\}$. By the pathwise uniqueness to \eqref{SMCH2 mollify}, we conclude that $z^R(t)=z^{\bar{R}}(t), \ \ t\in[0,\tau^R\wedge \tau^{\bar{R}})$. In particular, $\tau_R$ is increasing in $R$. Let $\tau^\ast=\lim_{R\to\infty}\tau_R$ and define
$$z=\sum_{R=1}^\infty1_{[\tau_{R-1},\tau_R)}z^R.$$
Then $(z,\tau^\ast)$ is the unique pathwise solution to \eqref{SMCH21} for $s>\frac{5}{2}$.

Next, we extend the range of $s$ to $s>\frac{3}{2}$. When $z_{0}\in L^{\infty}(\Omega, H^{s}\times H^{s})$ with $s>3/2$, by mollifying the initial data, we obtain a sequence of regular solutions $\{z_{n},\zeta^n\}_{n\in\mathbb N}$ to \eqref{SMCH21}. Motivated by \cite{GZ}, one can prove that there is some stopping time $\tau$ with $\mathbb{P}(\tau>0)=1$, a subsequence (still denoted by $z_{n}$) and some process $z$ such that $\mathbb P-a.s.$
\begin{equation*}\label{un}
\lim_{n\rightarrow\infty}\sup_{t\in[0,s]}(\|z_{n}-z\|_{H^{s}})=0~,\ \ s<\tau
\end{equation*}
and
\begin{equation}\label{u}
\sup_{t\in[0,s]}\|z\|_{H^{s}\times H^{s}}\leq\|z_{0}\|_{H^{s}\times H^{s}}+2, \ \ s<\tau.
\end{equation}
Then we can let $n\rightarrow\infty$ to prove that $(z,\tau)$ is a solution to \eqref{SMCH21}.

Besides, a cutting argument as in \cite{NV,DEM,GZ} enables us to remove the $L^{\infty}(\Omega, H^{s}\times H^{s})$ assumption on $(u_{0}, \gamma_{0})$. More precisely, consider the decomposition
\begin{equation*}\label{Omega}
\Omega_{m}=\{m-1\leq\|u_{0}\|_{H^{s}}+\|\gamma_{0}\|_{H^{s}}<m\},~~m\geq1.
\end{equation*}
We conclude $\sum_{m=1}^\infty \mathbb{P}(\Omega_{m})=1$.
Therefore we have $\mathbb P-a.s.$
\begin{equation*}\label{u0}
z_{0}(\omega,x)=\sum_{m\geq1}z^m_{0}(\omega,x):=\sum_{m\geq1}z_{0}(\omega,x)1_{\Omega_{m}}.
\end{equation*}
For each initial value $z^m_{0}$, we let $(z^{m}, \zeta_{m})$ be the pathwise unique solution to \eqref{SMCH21} satisfying \eqref{u}. Moreover, as $\Omega_{m}\cap\Omega_{m^{'}}=\emptyset, m\neq m^{'}$, $F_{1}(0,0)=0, F_{2}(0,0)=0$ and $h_{1}(t,0,0)=0, h_{2}(t,0,0)=0$ (see \eqref{h1}), it follows that
\begin{equation*}\label{u1}
z:=\sum_{m\geq1}z^m1_{\Omega_m}, \ \ \zeta=\sum_{m\geq1}\zeta_{m}1_{\Omega_m}
\end{equation*}
is the unique pathwise solution to \eqref{SMCH21} with corresponding initial condition $z_{0}$. Since $(z^{m}, \zeta_{m})$ satisfies \eqref{u}, we have $\mathbb{P}-a.s.$
\begin{eqnarray*}\label{u2}
\sup_{t\in[0,s]}(\|u\|_{H^{s}}^{2}+\|\gamma\|_{H^{s}}^{2})&=&\sum_{m=1}^{\infty}1_{\Omega_m} \sup_{t\in[0,s]}(\|u^{m}\|_{H^{s}}^{2}+\|\gamma^{m}\|_{H^{s}}^{2})\nonumber\\
&\leq&C\sum_{m=1}^{\infty}1_{\Omega_m} (4+\|u^m_{0}\|_{H^{s}}^{2}+ \|\gamma_{0}^m\|_{H^{s}}^{2})\nonumber\\
&=&C(4+\|u_{0}\|_{H^{s}}^{2}+\|\gamma_{0}\|_{H^{s}}^{2}),\ \ s<\zeta.
\end{eqnarray*}
So, \eqref{z} holds. Since the passage from $(z,\zeta)$ to a unique maximal pathwise solution $(z,\tau^{\ast})$ in the sense of Definition \ref{Pathwise} can be carried out as in \cite{NV,GZ,RZZ}, we omit the details. To finish the proof of Theorem \ref{Maximal solutions}, we only need to prove the blow-up scenario \eqref{Blow-up criterion}. Motivated by \cite{CH,DFD}, we consider the relationship between the explosion time of $\|z(t)\|_{H^{s}\times H^{s}}$ and the explosion time of $\|z(t)\|_{W^{1,\infty}\times W^{1,\infty}}$ in the next lemma.

\begin{lem}\label{blp}(Blow-up scenario 1)
Let $(z, \tau^{\ast})$ be the unique maximal solution to \eqref{SMCH21}. Then the real-valued stochastic processes $\|z(t)\|_{W^{1,\infty}\times W^{1,\infty}},~\|z(t)\|_{H^{s}\times H^{s}}$ are also $\mathcal{F}_{t}$-adapted. Besides, for any $m,n\in\mathbb Z^{+}$, define
\begin{eqnarray*}
\tau_{1,m}=\inf\{t\geq0: \|u(t)\|_{H^{s}}+\|\gamma(t)\|_{H^{s}}\geq m\}, ~\tau_{2,n}=\inf\{t\geq0: \|u(t)\|_{W^{1,\infty}}+\|\gamma(t)\|_{W^{1,\infty}}\geq n\}.
\end{eqnarray*}
For $\tau_{1}:=\tau^{\ast}=\lim_{m\rightarrow\infty}\tau_{1,m}$ and $\tau_{2}=\lim_{n\rightarrow\infty}\tau_{2,n}$, we have
\begin{eqnarray*}
\tau_{1}=\tau_{2}~~~\mathbb P-a.s.
\end{eqnarray*}
\end{lem}
Consequently, $1_{\{\lim_{t\rightarrow\tau^{\ast}}\|u(t)\|_{W^{1,\infty}}+\|\gamma(t)\|_{W^{1,\infty}}=\infty\}}=1_{\{\tau^{\ast}<\infty\}}~~\mathbb P-a.s.$.

\begin{pf}
Since $z\in C([0,\tau^\ast); H^{s}\times H^{s})$ almost surely, by the continuous embedding $H^{s}\times H^{s}\hookrightarrow W^{1,\infty}\times W^{1,\infty}$ for $s>3/2$, we conclude that $\|z(t)\|_{W^{1,\infty}\times W^{1,\infty}}$ is $\mathcal{F}_{t}$-adapted.  Moreover, the embedding $H^{s}\times H^{s}\hookrightarrow W^{1,\infty}\times W^{1,\infty}$ for $s>3/2$ means that $\mathbb P$-a.s. $\tau_{1}\leq\tau_{2}~\mathbb P-a.s.$ Now we only need to prove $\mathbb P$-a.s. $\tau_{2}\leq\tau_{1}$.
We cannot directly apply the It\^{o} formula for $\|u(t)\|_{H^{s}}^{2}+\|\gamma(t)\|_{H^{s}}^{2}$ since we only have $u, \gamma\in H^{s}$ and $uu_{x}, u\gamma_{x}\in H^{s-1}$. Therefore, the It\^{o} formula in Hilbert space cannot be applied directly, see (\cite{DG}, Theorem 4.32) or (\cite{GM}, Theorem 2.10). Instead, we will use the mollifier operator $T_{\epsilon}$ defined in \eqref{T} to overcome this difficulty. We apply $T_{\epsilon}$ to \eqref{SMCH21}, and then use the It\^{o}  formula for $\|T_{\epsilon}u\|_{H^{s}}^{2}, \|T_{\epsilon}\gamma\|_{H^{s}}^{2}$ to derive
\begin{eqnarray*}
d\|T_{\epsilon}u(t)\|_{H^{s}}^{2}&=&2(T_{\epsilon}h_{1}(u,\gamma)d\mathcal{W}_{1}, T_{\epsilon}u)_{H^{s}}-2(D^{s}T_{\epsilon}[uu_{x}], D^{s}T_{\epsilon}u)_{L^{2}}dt-2(D^{s}T_{\epsilon}F_{1}(u,\gamma), D^{s}T_{\epsilon}u)_{L^{2}}dt\nonumber\\
&+&\|T_{\epsilon}h_{1}(u,\gamma)\|^{2}_{L_2(U,H^{s})}dt,\nonumber\\
d\|T_{\epsilon}\gamma(t)\|_{H^{s}}^{2}&=&2(T_{\epsilon}h_{2}(u,\gamma)d\mathcal{W}_{2}, T_{\epsilon}\gamma)_{H^{s}}-2(D^{s}T_{\epsilon}[u\gamma_{x}], D^{s}T_{\epsilon}\gamma)_{L^{2}}dt-2(D^{s}T_{\epsilon}F_{2}(u,\gamma), D^{s}T_{\epsilon}\gamma)_{L^{2}}dt\nonumber\\
&+&\|T_{\epsilon}h_{2}(u,\gamma)\|^{2}_{L_2(U,H^{s})}dt.
\end{eqnarray*}
Therefore for any $n_{1},  m \geq 1$, $r\geq 0$, and $t\in [0, \tau_{2,n_{1}}\wedge r\wedge \tau_{1,m}]$,
\begin{align*}
&\|T_{\epsilon}u(t)\|_{H^{s}}^{2}+\|T_{\epsilon}\gamma(t)\|_{H^{s}}^{2}-\|T_{\epsilon}u(0)\|_{H^{s}}^{2}-\|T_{\epsilon}\gamma(0)\|_{H^{s}}^{2}\nonumber\\
=&2\sum_{j=1}^{\infty}\int_{0}^{t}(D^{s}T_{\epsilon}h_{1}(u,\gamma)e_{j}, D^{s}T_{\epsilon}u)_{L^{2}}dW_{j}^{1}+2\sum_{i=1}^{\infty}\int_{0}^{t}(D^{s}T_{\epsilon}h_{2}(u,\gamma)e_{i}, D^{s}T_{\epsilon}\gamma)_{L^{2}}dW_{i}^{2}\nonumber\\
&-2\int_{0}^{t}(D^{s}T_{\epsilon}[uu_{x}], D^{s}T_{\epsilon}u)_{L^{2}}dt^{'}-2\int_{0}^{t}(D^{s}T_{\epsilon}F_{1}(u,\gamma), D^{s}T_{\epsilon}u)_{L^{2}}dt^{'}\nonumber\\
&+\int_{0}^{t}\|T_{\epsilon}h_{1}(u,\gamma)\|^{2}_{L_2(U, H^{s})}dt^{'}
-2\int_{0}^{t}(D^{s}T_{\epsilon}[u\gamma_{x}], D^{s}T_{\epsilon}\gamma)_{L^{2}}dt^{'}\nonumber\\
&-2\int_{0}^{t}(D^{s}T_{\epsilon}F_{2}(u,\gamma), D^{s}T_{\epsilon}\gamma)_{L^{2}}dt^{'}+\int_{0}^{t}\|T_{\epsilon}h_{2}(u,\gamma)\|^{2}_{L_2(U, H^{s})}dt^{'}\nonumber\\
=:&\int_{0}^{t}\sum_{j=1}^{\infty}L_{1,j}dW_{j}^{1}+ \int_{0}^{t}\sum_{i=1}^{\infty}L_{2,i}dW_{i}^{2}+\sum_{j=3}^{8}\int_{0}^{t}L_{j}dt^{'},
\end{align*}
where $\{e_{k}\}$ is the complete orthonormal basis of $U$. On account of the Burkholder-Davis-Gundy inequality and  \eqref{h1}, we obtain that
\begin{flalign*}
&\mathbb E\bigg[\sup_{t\in[0, \tau_{2,n_{1}}\wedge r\wedge \tau_{1,m}]}\left|\int_{0}^{t}\sum_{j=1}^{\infty}L_{1,j}dW_{j}^{1}\right|\bigg|\mathcal{F}_0\bigg]\leq C\mathbb E\bigg\{\bigg(\sum_{j=1}^{\infty}\int_{0}^{\tau_{2,n_{1}}\wedge r\wedge \tau_{1,m}}|L_{1,j}|^{2}d t\bigg)^{\frac{1}{2}}\bigg|\mathcal{F}_0\bigg\} \nonumber\\
&\leq\frac{1}{2}\mathbb E\left(\sup_{t\in[0, \tau_{2,n_{1}}\wedge r\wedge \wedge \tau_{1,m}]}\|T_{\epsilon}u\|^{2}_{H^{s}}\bigg|\mathcal{F}_0\right)+Cf^{2}(n_{1})\int_{0}^{r}\mathbb E\bigg[\sup_{t^{'}\in[0, \tau_{2,n_{1}}\wedge t\wedge \tau_{1,m}]}(1+\|u(t^{'})\|^{2}_{H^{s}}+\|\gamma(t^{'})\|^{2}_{H^{s}})\bigg|\mathcal{F}_0\bigg]d t,  \nonumber\\
&\mathbb E\bigg[\sup_{t\in[0, \tau_{2,n_{1}}\wedge r\wedge \tau_{1,m}]}\left|\int_{0}^{t}\sum_{i=1}^{\infty}L_{2,i}dW_{i}^{2}\right|\bigg|\mathcal{F}_0\bigg] \leq C\mathbb E\bigg\{\bigg(\sum_{i=1}^{\infty}\int_{0}^{\tau_{2,n_{1}}\wedge r\wedge \tau_{1,m}}|L_{1,i}|^{2}d t\bigg)^{\frac{1}{2}}\bigg|\mathcal{F}_0\bigg\}\nonumber\\
&\leq\frac{1}{2}\mathbb E\left[\sup_{t\in[0, \tau_{2,n_{1}}\wedge r\wedge \tau_{1,m}]}\|T_{\epsilon}\gamma\|^{2}_{H^{s}}\bigg|\mathcal{F}_0\right]+Cf^{2}(n_{1})\int_{0}^{r}\mathbb{E}\bigg[\sup_{t^{'}\in[0, \tau_{2,n_{1}}\wedge t\wedge \tau_{1,m}]}(1+\|u(t^{'})\|^{2}_{H^{s}}+\|\gamma(t^{'})\|^{2}_{H^{s}})\bigg|\mathcal{F}_0\bigg]d t.
\end{flalign*}
For $L_{3}, L_{6}$, using integration by part, Sobolev's inequality and Lemma \ref{l2}, we have
\begin{eqnarray*}
(D^{s}T_{\epsilon}[uu_{x}], D^{s}T_{\epsilon}u)_{L^{2}}&=&([D^{s},u]u_{x}, D^{s}T^{2}_{\epsilon}u)_{L^{2}}+([T_{\epsilon},u]D^{s}u_{x}, D^{s}T_{\epsilon}u)_{L^{2}}+(uD^{s}T_{\epsilon}u_{x}, D^{s}T_{\epsilon}u)_{L^{2}}\nonumber\\
&\leq&C\|u\|_{W^{1,\infty}}\|u\|^{2}_{H^{s}},\nonumber\\
(D^{s}T_{\epsilon}[u\gamma_{x}], D^{s}T_{\epsilon}\gamma)_{L^{2}}&=&([D^{s},u]\gamma_{x}, D^{s}T^{2}_{\epsilon}\gamma)_{L^{2}}+([T_{\epsilon},u]D^{s}\gamma_{x}, D^{s}T_{\epsilon}\gamma)_{L^{2}}+(uD^{s}T_{\epsilon}\gamma_{x}, D^{s}T_{\epsilon}\gamma)_{L^{2}}\nonumber\\
&\leq&C(\|u\|_{W^{1,\infty}}+\|\gamma\|_{W^{1,\infty}})(\|u\|^{2}_{H^{s}}+\|\gamma\|^{2}_{H^{s}}).
\end{eqnarray*}
For $L_{4}, L_{7}$, we derive from Lemma \ref{F1,F2} that
\begin{eqnarray*}
(D^{s}T_{\epsilon}F_{1}(u,\gamma), D^{s}T_{\epsilon}u)_{L^{2}}&\leq&C(\|u\|_{W^{1,\infty}}+\|\gamma\|_{W^{1,\infty}})(\|u\|^{2}_{H^{s}}+\|\gamma\|^{2}_{H^{s}}),\nonumber\\
(D^{s}T_{\epsilon}F_{2}(u,\gamma), D^{s}T_{\epsilon}\gamma)_{L^{2}}&\leq&C(\|u\|_{W^{1,\infty}}+\|\gamma\|_{W^{1,\infty}})(\|u\|^{2}_{H^{s}}+\|\gamma\|^{2}_{H^{s}}).
\end{eqnarray*}
For $L_{5}, L_{8}$, it follows from the Assumption \ref{h1} that
\begin{align*}
&\mathbb E\bigg\{\int_{0}^{\tau_{2,n_{1}}\wedge r\wedge \tau_{1,m}}\|T_{\epsilon}h_{1}(u,\gamma)\|^{2}_{L_2(U, H^{s})}dt^{'}\bigg|\mathcal{F}_0\bigg\}\\
&\leq Cf^{2}(n_{1})\int_{0}^{r}\mathbb{E}\bigg\{\sup_{t^{'}\in[0, \tau_{2,n_{1}}\wedge t\wedge \tau_{1,m}]}(1+\|u(t^{'})\|^{2}_{H^{s}}+\|\gamma(t^{'})\|^{2}_{H^{s}})\bigg|\mathcal{F}_0\bigg\}d t,\\
&\mathbb E\bigg\{\int_{0}^{\tau_{2,n_{1}}\wedge r\wedge \tau_{1,m}}\|T_{\epsilon}h_{2}(u,\gamma)\|^{2}_{L_2(U, H^{s})}dt^{'}\bigg|\mathcal{F}_0\bigg\}\\
&\leq Cf^{2}(n_{1})\int_{0}^{r}\mathbb{E}\bigg\{\sup_{t^{'}\in[0, \tau_{2,n_{1}}\wedge t\wedge \tau_{1,m}]}(1+\|u(t^{'})\|^{2}_{H^{s}}+\|\gamma(t^{'})\|^{2}_{H^{s}})\bigg|\mathcal{F}_0\bigg\}d t.
\end{align*}
Therefore combining the above estimates, we have
\begin{eqnarray*}
&&\mathbb E\left\{\sup_{t\in[0, \tau_{2,n_{1}}\wedge r\wedge \tau_{1,m}]}(\|T_{\epsilon}u(t)\|_{H^{s}}^{2}+\|T_{\epsilon}\gamma(t)\|_{H^{s}}^{2})\bigg|\mathcal{F}_0\right\}\\
&\leq& C [\|u(0)\|_{H^{s}}^{2}+\|\gamma(0)\|_{H^{s}}^{2}]+C\int_{0}^{r}\mathbb E\bigg(1+\sup_{t^{'}\in[0,\tau_{2,n_{1}}\wedge t\wedge \tau_{1,m}]}(\|u(t^{'})\|^{2}_{H^{s}}+\|\gamma(t^{'})\|^{2}_{H^{s}})\bigg|\mathcal{F}_0\bigg)d t.
\end{eqnarray*}
Since the right hand side of the above estimate does not depend on $\epsilon$, and $(T_{\epsilon}u, T_{\epsilon}\gamma)$ tends to $(u, \gamma)$ in $C([0, \tau_{1,m}\wedge r], H^{s}\times H^{s})$ almost surely as $\epsilon\rightarrow0$, by Fatou's Lemma, one can send $\epsilon\rightarrow0$ to obtain
\begin{align*}
&\mathbb E\left\{\sup_{t\in[0, \tau_{2,n_{1}}\wedge r\wedge \tau_{1,m}]}(\|u(t)\|_{H^{s}}^{2}+\|\gamma(t)\|_{H^{s}}^{2})\bigg|\mathcal{F}_0\right\}\nonumber \\
&\leq C [\|u(0)\|_{H^{s}}^{2}+\|\gamma(0)\|_{H^{s}}^{2}]+C\int_{0}^{r}\mathbb E\bigg(1+\sup_{t^{'}\in[0,\tau_{2,n_{1}}\wedge t\wedge \tau_{1,m}]}(\|u(t^{'})\|^{2}_{H^{s}}+\|\gamma(t^{'})\|^{2}_{H^{s}})|\mathcal{F}_0\bigg)d t.
\end{align*}
Then Gronwall's inequality shows that for each $n_{1}\in \mathbb Z^{+}$, $r\in\mathbb R^{+}$, there is a constant $C=C(n_{1}, r, u_{0}, \gamma_{0})>0$ such that
\begin{eqnarray*}
\mathbb E\left[\sup_{t\in[0, \tau_{2,n_{1}}\wedge r\wedge \tau_{1,m}]}(\|u(t)\|_{H^{s}}^{2}+\|\gamma(t)\|_{H^{s}}^{2})\bigg|\mathcal{F}_0\right]<C(n_{1}, r, u_{0}, \gamma_{0}).
\end{eqnarray*}
So, it follows from Chebyshev's inequality and Fatou's lemma that
\begin{align*}\mathbb{P}(\tau_{1}\leq \tau_{2,n_{1}}\wedge r|\mathcal{F}_0)&\leq \lim_{m\to\infty}\mathbb{P}(\tau_{1,m}
\leq \tau_{2,n_{1}}\wedge r|\mathcal{F}_0)\\
&\leq \lim_{m\to\infty}\mathbb{P}\left(\sup_{t\in[0, \tau_{2,n_{1}}\wedge r\wedge \tau_{1,m}]}(\|u(t)\|_{H^{s}}+\|\gamma(t)\|_{H^{s}})\geq m\bigg|\mathcal{F}_0\right)\\
&\leq \lim_{m\to\infty}\frac{\mathbb E[\sup_{t\in[0, \tau_{2,n_{1}}\wedge r\wedge \tau_{1,m}]}(2\|u(t)\|_{H^{s}}^{2}+2\|\gamma(t)\|_{H^{s}}^{2})|\mathcal{F}_0]}{m^2}=0.
\end{align*}
Letting $r\to\infty$ first and then $n_1\to\infty$, Fatou's lemma yields that
$$\mathbb{P}(\tau_{1}\leq \tau_{2}|\mathcal{F}_0)=0.$$
Therefore, we conclude
$$\mathbb{P}(\tau_{1}\geq  \tau_{2})\geq 1-\mathbb{P}(\tau_{1}\leq \tau_{2})=1-\mathbb{P}[\mathbb{P}(\tau_{1}\leq \tau_{2}|\mathcal{F}_0)]=1.$$
We finish the section with the proof of the first blow-up scenario \eqref{Blow-up criterion}.
\end{pf}

\section{Proof of Theorem \ref{Global strong nonlinear}}
Assume $s>5/2$ and let $(u_{0}, \gamma_{0})$ be $H^{s}\times H^{s}$-valued $\mathcal{F}_{0}$-measurable random variable. Let $h_{1}(t, u, \gamma)=a(t)(1+\|u\|_{W^{1,\infty}}+\|\gamma\|_{W^{1,\infty}})^{\theta}u, h_{2}(t, u, \gamma)=a(t)(1+\|u\|_{W^{1,\infty}}+\|\gamma\|_{W^{1,\infty}})^{\theta}\gamma$ with $\theta\geq1/2$ and $a(t)\neq0$.

For $s>3/2$, the embedding $H^{s}\times H^{s}\hookrightarrow W^{1,\infty}\times W^{1,\infty}$ implies that
\begin{eqnarray*}\label{h22}
\sup_{\|u_{1}\|_{H^{s}},\|\gamma_{1}\|_{H^{s}},\|u_{2}\|_{H^{s}},\|\gamma_{2}\|_{H^{s}}\leq N}&&\sum_{i=1}^2(\|h_{i}(t, u_{1},\gamma_{1})-h_{i}(t, u_{2},\gamma_{2})\|_{H^{s}}\nonumber\\
&\leq& g(N)(\|u_{1}-u_{2}\|_{H^{s}}+\|\gamma_{1}-\gamma_{2}\|_{H^{s}}),~~N\geq1.
\end{eqnarray*}
 Hence, by Theorem \ref{Maximal solutions}, we conclude that \eqref{SMCH21} admits a unique pathwise solution $z=(u,\gamma)$ in $H^{s}\times H^{s}$ with $s>5/2$ and maximal existence time $\tau^{\ast}$. Define
\begin{eqnarray*}\label{time}
\tau_{m}=\inf\{t\geq0: \|u\|_{H^{s}}^{2}+\|\gamma\|_{H^{s}}^{2}\geq m\}.
\end{eqnarray*}
Applying the It\^{o} formula to $\|T_{\epsilon}u\|_{H^{s}}^{2}, \|T_{\epsilon}\gamma\|_{H^{s}}^{2}$ gives
\begin{align*}
d\|T_{\epsilon}u\|_{H^{s}}^{2}=&2a(t)(1+\|u\|_{W^{1,\infty}}+\|\gamma\|_{W^{1,\infty}})^{\theta}(T_{\epsilon}u,T_{\epsilon}u)_{H^{s}}d W-2(T_{\epsilon}[uu_{x}],T_{\epsilon}u)_{H^{s}}d t\nonumber\\
&-2(T_{\epsilon}F_{1}(u,\gamma),T_{\epsilon}u)_{H^{s}}d t
+a^{2}(t)(1+\|u\|_{W^{1,\infty}}+\|\gamma\|_{W^{1,\infty}})^{2\theta}(T_{\epsilon}u,T_{\epsilon}u)_{H^{s}}d t,\nonumber\\
d\|T_{\epsilon}\gamma\|_{H^{s}}^{2}=&2a(t)(1+\|u\|_{W^{1,\infty}}+\|\gamma\|_{W^{1,\infty}})^{\theta}(T_{\epsilon}\gamma,T_{\epsilon}\gamma)_{H^{s}}d W-2(T_{\epsilon}[u\gamma_{x}],T_{\epsilon}\gamma)_{H^{s}}d t\nonumber\\
&-2(T_{\epsilon}F_{2}(u,\gamma),T_{\epsilon}\gamma)_{H^{s}}d t
+a^{2}(t)(1+\|u\|_{W^{1,\infty}}+\|\gamma\|_{W^{1,\infty}})^{2\theta}(T_{\epsilon}\gamma,T_{\epsilon}\gamma)_{H^{s}}d t,
\end{align*}
Again, using It\^o's formula to $\log(1+\|T_{\epsilon}u\|_{H^{s}}^{2}+\|T_{\epsilon}\gamma\|_{H^{s}}^{2})$ yields
\begin{align*}
&d \log(1+\|T_{\epsilon}u\|_{H^{s}}^{2}+\|T_{\epsilon}\gamma\|_{H^{s}}^{2})
=\frac{2a(t)(1+\|u\|_{W^{1,\infty}}+\|\gamma\|_{W^{1,\infty}})^{\theta}(T_{\epsilon}u,T_{\epsilon}u)_{H^{s}}}{1+\|T_{\epsilon}u\|_{H^{s}}^{2}+\|T_{\epsilon}\gamma\|_{H^{s}}^{2}}d W \nonumber\\
&-\frac{2(T_{\epsilon}[uu_{x}],T_{\epsilon}u)_{H^{s}}}{1+\|T_{\epsilon}u\|_{H^{s}}^{2}+\|T_{\epsilon}\gamma\|_{H^{s}}^{2}}d t-\frac{2(T_{\epsilon}F_{1}(u,\gamma),T_{\epsilon}u)_{H^{s}}}{1+\|T_{\epsilon}u\|_{H^{s}}^{2}+\|T_{\epsilon}\gamma\|_{H^{s}}^{2}}d t \nonumber\\
&+\frac{a^{2}(t)(1+\|u\|_{W^{1,\infty}}+\|\gamma\|_{W^{1,\infty}})^{2\theta}(T_{\epsilon}u,T_{\epsilon}u)_{H^{s}}}{1+\|T_{\epsilon}u\|_{H^{s}}^{2}+\|T_{\epsilon}\gamma\|_{H^{s}}^{2}}d t\nonumber\\
&+\frac{2a(t)(1+\|u\|_{W^{1,\infty}}+\|\gamma\|_{W^{1,\infty}})^{\theta}(T_{\epsilon}\gamma,T_{\epsilon}\gamma)_{H^{s}}}{1+\|T_{\epsilon}u\|_{H^{s}}^{2}+\|T_{\epsilon}\gamma\|_{H^{s}}^{2}}d W \nonumber\\
&-\frac{2(T_{\epsilon}[u\gamma_{x}],T_{\epsilon}\gamma)_{H^{s}}}{1+\|T_{\epsilon}u\|_{H^{s}}^{2}+\|T_{\epsilon}\gamma\|_{H^{s}}^{2}}d t-\frac{2(T_{\epsilon}F_{2}(u,\gamma),T_{\epsilon}\gamma)_{H^{s}}}{1+\|T_{\epsilon}u\|_{H^{s}}^{2}+\|T_{\epsilon}\gamma\|_{H^{s}}^{2}}d t\nonumber\\
&+\frac{a^{2}(t)(1+\|u\|_{W^{1,\infty}}+\|\gamma\|_{W^{1,\infty}})^{2\theta}(T_{\epsilon}\gamma,T_{\epsilon}\gamma)_{H^{s}}}{1+\|T_{\epsilon}u\|_{H^{s}}^{2}+\|T_{\epsilon}\gamma\|_{H^{s}}^{2}}d t\nonumber\\
&-2a^{2}(t)(1+\|u\|_{W^{1,\infty}}+\|\gamma\|_{W^{1,\infty}})^{2\theta}\nonumber\\
&\qquad\qquad\qquad\qquad\quad\times \frac{[(T_{\epsilon}u,T_{\epsilon}u)_{H^{s}}^{2}+(T_{\epsilon}\gamma,T_{\epsilon}\gamma)_{H^{s}}^{2}+2(T_{\epsilon}u,T_{\epsilon}u)_{H^{s}}(T_{\epsilon}\gamma,T_{\epsilon}\gamma)_{H^{s}}]}{(1+\|T_{\epsilon}u\|_{H^{s}}^{2}+\|T_{\epsilon}\gamma\|_{H^{s}}^{2})^{2}}d t.
\end{align*}
By Lemma \ref{Te}-\ref{F1,F2} and Lemma \ref{energy estimation}, it follows that
\begin{align*}
&\mathbb E[\log(1+\|T_{\epsilon}u(t\wedge\tau_{m})\|_{H^{s}}^{2}+\|T_{\epsilon}\gamma(t\wedge\tau_{m})\|_{H^{s}}^{2})|\mathcal{F}_0]
-\log(1+\|T_{\epsilon}u_{0}\|_{H^{s}}^{2}+\|T_{\epsilon}\gamma_{0}\|_{H^{s}}^{2})\nonumber\\
=&-2\mathbb E\bigg[\int_{0}^{t\wedge\tau_{m}}\frac{(T_{\epsilon}[uu_{x}],T_{\epsilon}u)_{H^{s}}}{1+\|T_{\epsilon}u\|_{H^{s}}^{2}+\|T_{\epsilon}\gamma\|_{H^{s}}^{2}}d t^{'}\bigg|\mathcal{F}_0\bigg]-2\mathbb E\bigg[\int_{0}^{t\wedge\tau_{m}}\frac{(T_{\epsilon}F_{1}(u,\gamma),T_{\epsilon}u)_{H^{s}}}{1+\|T_{\epsilon}u\|_{H^{s}}^{2}+\|T_{\epsilon}\gamma\|_{H^{s}}^{2}}d t^{'}\bigg|\mathcal{F}_0\bigg]\nonumber\\
&+\mathbb E\bigg[\int_{0}^{t\wedge\tau_{m}}\frac{a^{2}(t^{'})(1+\|u\|_{W^{1,\infty}}+\|\gamma\|_{W^{1,\infty}})^{2\theta}(T_{\epsilon}u,T_{\epsilon}u)_{H^{s}}}{1+\|T_{\epsilon}u\|_{H^{s}}^{2}+\|T_{\epsilon}\gamma\|_{H^{s}}^{2}}d t^{'}\bigg|\mathcal{F}_0\bigg]\nonumber\\
&-2\mathbb E\bigg[\int_{0}^{t\wedge\tau_{m}}\frac{(T_{\epsilon}[u\gamma_{x}],T_{\epsilon}\gamma)_{H^{s}}}{1+\|T_{\epsilon}u\|_{H^{s}}^{2}+\|T_{\epsilon}\gamma\|_{H^{s}}^{2}}d t^{'}\bigg|\mathcal{F}_0\bigg]-2\mathbb E\bigg[\int_{0}^{t\wedge\tau_{m}}\frac{(T_{\epsilon}F_{2}(u,\gamma),T_{\epsilon}\gamma)_{H^{s}}}{1+\|T_{\epsilon}u\|_{H^{s}}^{2}+\|T_{\epsilon}\gamma\|_{H^{s}}^{2}}d t^{'}\bigg|\mathcal{F}_0\bigg]\nonumber\\
&+\mathbb E\bigg[\int_{0}^{t\wedge\tau_{m}}\frac{a^{2}(t^{'})(1+\|u\|_{W^{1,\infty}}+\|\gamma\|_{W^{1,\infty}})^{2\theta}(T_{\epsilon}\gamma,T_{\epsilon}\gamma)_{H^{s}}}{1+\|T_{\epsilon}u\|_{H^{s}}^{2}+\|T_{\epsilon}\gamma\|_{H^{s}}^{2}}d t^{'}\bigg|\mathcal{F}_0\bigg]\nonumber\\
&-2\mathbb E\bigg[\int_{0}^{t\wedge\tau_{m}}a^{2}(t^{'})(1+\|u\|_{W^{1,\infty}}+\|\gamma\|_{W^{1,\infty}})^{2\theta}\nonumber\\
&\qquad\qquad\qquad\qquad\quad\times \frac{[(T_{\epsilon}u,T_{\epsilon}u)_{H^{s}}^{2}+(T_{\epsilon}\gamma,T_{\epsilon}\gamma)_{H^{s}}^{2}++2(T_{\epsilon}u,T_{\epsilon}u)_{H^{s}}(T_{\epsilon}\gamma,T_{\epsilon}\gamma)_{H^{s}}]}{(1+\|T_{\epsilon}u\|_{H^{s}}^{2}+\|T_{\epsilon}\gamma\|_{H^{s}}^{2})^{2}}d t^{'}\bigg|\mathcal{F}_0\bigg]\nonumber\\
\leq&\mathbb E\bigg[\int_{0}^{t\wedge\tau_{m}}\frac{K(\|u\|_{W^{1,\infty}}+\|\gamma\|_{W^{1,\infty}})(\|u\|_{H^{s}}^{2}+\|\gamma\|_{H^{s}}^{2})}{1+\|T_{\epsilon}u\|_{H^{s}}^{2}+\|T_{\epsilon}\gamma\|_{H^{s}}^{2}}d t^{'}\bigg|\mathcal{F}_0\bigg]\nonumber\\
&+\mathbb E\bigg[\int_{0}^{t\wedge\tau_{m}}\frac{a^{2}(t^{'})(1+\|u\|_{W^{1,\infty}}+\|\gamma\|_{W^{1,\infty}})^{2\theta}(\|T_{\epsilon}u\|_{H^{s}}^{2}+\|T_{\epsilon}\gamma\|_{H^{s}}^{2})}{1+\|T_{\epsilon}u\|_{H^{s}}^{2}+\|T_{\epsilon}\gamma\|_{H^{s}}^{2}}d t^{'}\bigg|\mathcal{F}_0\bigg]\nonumber\\
&-2\mathbb E\bigg[\int_{0}^{t\wedge\tau_{m}}a^{2}(t^{'})(1+\|u\|_{W^{1,\infty}}+\|\gamma\|_{W^{1,\infty}})^{2\theta}\nonumber\\
&\qquad\qquad\qquad\qquad\quad\frac{(\|T_{\epsilon}u\|_{H^{s}}^{4}+\|T_{\epsilon}\gamma\|_{H^{s}}^{4}+2\|T_{\epsilon}u\|_{H^{s}}^{2}\|T_{\epsilon}\gamma\|_{H^{s}}^{2})}{(1+\|T_{\epsilon}u\|_{H^{s}}^{2}+\|T_{\epsilon}\gamma\|_{H^{s}}^{2})^{2}}d t^{'}\bigg|\mathcal{F}_0\bigg].
\end{align*}
 Let
\begin{align}\label{dom}
&I^\epsilon_1(t')=\frac{K(\|u\|_{W^{1,\infty}}+\|\gamma\|_{W^{1,\infty}})(\|u\|_{H^{s}}^{2}+ \|\gamma\|_{H^{s}}^{2})}{1+\|T_{\epsilon}u\|_{H^{s}}^{2}+\|T_{\epsilon}\gamma\|_{H^{s}}^{2}} -\frac{K(\|u\|_{W^{1,\infty}}+\|\gamma\|_{W^{1,\infty}}) (\|u\|_{H^{s}}^{2}+\|\gamma\|_{H^{s}}^{2})}{1+\|u\|_{H^{s}}^{2}+\|\gamma\|_{H^{s}}^{2}} \nonumber,\\
&I^\epsilon_2(t')=\frac{a^{2}(t^{'})(1+\|u\|_{W^{1,\infty}}+\|\gamma\|_{W^{1,\infty}})^{2\theta} (\|T_{\epsilon}u\|_{H^{s}}^{2}+\|T_{\epsilon}\gamma\|_{H^{s}}^{2})}{1+\|T_{\epsilon}u\|_{H^{s}}^{2}+ \|T_{\epsilon}\gamma\|_{H^{s}}^{2}}\nonumber\\
&\qquad\qquad\qquad-\frac{a^{2}(t^{'})(1+\|u\|_{W^{1,\infty}}+\|\gamma\|_{W^{1,\infty}})^{2\theta}(\|u\|_{H^{s}}^{2}+ \|\gamma\|_{H^{s}}^{2})}{1+\|u\|_{H^{s}}^{2}+\|\gamma\|_{H^{s}}^{2}}\nonumber,\\
&I^\epsilon_3(t')=\frac{a^{2}(t^{'})(1+\|u\|_{W^{1,\infty}}+\|\gamma\|_{W^{1,\infty}})^{2\theta}(\|T_{\epsilon}u\|_{H^{s}}^{4}+\|T_{\epsilon}\gamma\|_{H^{s}}^{4}+2\|T_{\epsilon}u\|_{H^{s}}^{2}\|T_{\epsilon}\gamma\|_{H^{s}}^{2})}{(1+\|T_{\epsilon}u\|_{H^{s}}^{2}+\|T_{\epsilon}\gamma\|_{H^{s}}^{2})^{2}}\nonumber\\
&\qquad\qquad-\frac{a^{2}(t^{'})(1+\|u\|_{W^{1,\infty}}+\|\gamma\|_{W^{1,\infty}})^{2\theta}(\|u\|_{H^{s}}^{4}+\|\gamma\|_{H^{s}}^{4}+2\|u\|_{H^{s}}^{2}\|\gamma\|_{H^{s}}^{2})}{(1+\|u\|_{H^{s}}^{2}+\|\gamma\|_{H^{s}}^{2})^{2}}.
\end{align}
Notice that for any $T>0$, $(T_{\epsilon}u, T_{\epsilon}\gamma)$ tends to $(u, \gamma)$ in $C([0, \tau_m\wedge t], H^{s}\times H^{s})$ almost surely as $\epsilon\rightarrow0$. It follows from the dominated convergence theorem that
$$\lim_{\epsilon\to 0}\mathbb E\left[\int_{0}^{t\wedge\tau_{m}}[|I^\epsilon_1(t')|+|I_2^\epsilon(t')|+|I^\epsilon_3(t')|]dt'\bigg|\mathcal{F}_0\right]=0.$$
Then, by \eqref{T} and
the dominated convergence Theorem, it holds
\begin{align*}
&\mathbb E[\log(1+\|u(t\wedge\tau_{m})\|_{H^{s}}^{2}+\|\gamma(t\wedge\tau_{m})\|_{H^{s}}^{2})|\mathcal{F}_0]-\log(1+\|u_{0}\|_{H^{s}}^{2}+\|\gamma_{0}\|_{H^{s}}^{2})\nonumber\\
\leq & \mathbb E\bigg[\int_{0}^{t\wedge\tau_{m}}\frac{K(\|u\|_{W^{1,\infty}}+\|\gamma\|_{W^{1,\infty}})(\|u\|_{H^{s}}^{2}+\|\gamma\|_{H^{s}}^{2})}{1+\|u\|_{H^{s}}^{2}+\|\gamma\|_{H^{s}}^{2}}d t^{'}\bigg|\mathcal{F}_0\bigg]\nonumber\\
&+\mathbb E\bigg[\int_{0}^{t\wedge\tau_{m}}\frac{a^{2}(t^{'})(1+\|u\|_{W^{1,\infty}}+\|\gamma\|_{W^{1,\infty}})^{2\theta}(\|u\|_{H^{s}}^{2}+\|\gamma\|_{H^{s}}^{2})}{1+\|u\|_{H^{s}}^{2}+\|\gamma\|_{H^{s}}^{2}}d t^{'}\bigg|\mathcal{F}_0\bigg]\nonumber\\
&-2\mathbb E\bigg[\int_{0}^{t\wedge\tau_{m}}\frac{a^{2}(t^{'})(1+\|u\|_{W^{1,\infty}}+\|\gamma\|_{W^{1,\infty}})^{2\theta}(\|u\|_{H^{s}}^{4}+\|\gamma\|_{H^{s}}^{4}+2\|u\|_{H^{s}}^{2}\|\gamma\|_{H^{s}}^{2})}{(1+\|u\|_{H^{s}}^{2}+\|\gamma\|_{H^{s}}^{2})^{2}}d t^{'}\bigg|\mathcal{F}_0\bigg].
\end{align*}
Lemma \ref{lem strong noise} with  $x_{1}=\|u\|_{W^{1,\infty}}, x_{2}=\|\gamma\|_{W^{1,\infty}}, y_{1}=\|u\|_{H^{s}}, y_{2}=\|\gamma\|_{H^{s}}$ immediately shows that there are constants $C_{1}, C_{2}>0$ such
that
\begin{eqnarray*}
&&\mathbb E[\log(1+\|u(t\wedge\tau_{m})\|_{H^{s}}^{2}+\|\gamma(t\wedge\tau_{m})\|_{H^{s}}^{2})|\mathcal{F}_0]-\log(1+\|u_{0}\|_{H^{s}}^{2}+\|\gamma_{0}\|_{H^{s}}^{2})\nonumber\\
&\leq&\mathbb E\bigg[\int_{0}^{t\wedge\tau_{m}}C_{1}-C_{2}\frac{a^{2}(t^{'})(1+\|u\|_{W^{1,\infty}}+\|\gamma\|_{W^{1,\infty}})^{2\theta}(\|u\|_{H^{s}}^{2}+\|\gamma\|_{H^{s}}^{2})^{2}}{(1+\|u\|_{H^{s}}^{2}+\|\gamma\|_{H^{s}}^{2})^{2}(1+\log(1+\|u\|_{H^{s}}^{2}+\|\gamma\|_{H^{s}}^{2}))}d t^{'}\bigg|\mathcal{F}_0\bigg],
\end{eqnarray*}
which means that for some $C(u_{0}, \gamma_{0}, C_{1}, C_{2}, t)>0$,
\begin{eqnarray}\label{lemma3}
\mathbb E\bigg[\int_{0}^{t\wedge\tau_{m}}\frac{a^{2}(t^{'})(1+\|u\|_{W^{1,\infty}} +\|\gamma\|_{W^{1,\infty}})^{2\theta}(\|u\|_{H^{s}}^{2}+ \|\gamma\|_{H^{s}}^{2})^{2}}{(1+\|u\|_{H^{s}}^{2}+ \|\gamma\|_{H^{s}}^{2})^{2}(1+\log(1+\|u\|_{H^{s}}^{2}+\|\gamma\|_{H^{s}}^{2}))}d t^{'}\bigg|\mathcal{F}_0\bigg]\leq C(u_{0}, \gamma_{0}, C_{1}, C_{2}, t).
\end{eqnarray}
Therefore, for any $T>0$, by using Lemma \ref{lem strong noise} and the Burkholder-Davis-Gundy inequality, we find that
\begin{align*}
&\mathbb E\left[\sup_{t\in[0, T\wedge\tau_{m}]}\log(1+\|T_{\epsilon}u\|_{H^{s}}^{2}+\|T_{\epsilon}\gamma\|_{H^{s}}^{2}) \bigg|\mathcal{F}_0\right]-\log(1+\|T_{\epsilon}u_{0}\|_{H^{s}}^{2}+\|T_{\epsilon}\gamma_{0}\|_{H^{s}}^{2})\nonumber\\
\leq& C\left(\mathbb E\bigg[\int_{0}^{T\wedge\tau_{m}}\frac{a^{2}(t)(1+\|u\|_{W^{1,\infty}}+\|\gamma\|_{W^{1,\infty}})^{2\theta}(\|T_{\epsilon}u\|_{H^{s}}^{2}+\|T_{\epsilon}\gamma\|_{H^{s}}^{2})^{2}}{(1+\|T_{\epsilon}u\|_{H^{s}}^{2}+\|T_{\epsilon}\gamma\|_{H^{s}}^{2})^{2}}d t\bigg|\mathcal{F}_0\bigg]\right)^{\frac{1}{2}}\nonumber\\
&+\mathbb E\bigg[\int_{0}^{T\wedge\tau_{m}}\left|C_{1}-C_{2}\frac{a^{2}(t) (1+\|u\|_{W^{1,\infty}}+\|\gamma\|_{W^{1,\infty}})^{2\theta} (\|u\|_{H^{s}}^{2}+\|\gamma\|_{H^{s}}^{2})^{2}}{(1+\|u\|_{H^{s}}^{2}+ \|\gamma\|_{H^{s}}^{2})^{2}(1+\log(1+\|u\|_{H^{s}}^{2}+\|\gamma\|_{H^{s}}^{2}))}\right|d t\bigg|\mathcal{F}_0\bigg]\\
&+\mathbb{E}\left[\int_0^{T\wedge\tau_{m}}[|I^\epsilon_1(t)|+|I_2^\epsilon(t)|+|I^\epsilon_3(t)|]d t\bigg|\mathcal{F}_0\right]\nonumber\\
\leq&\frac{1}{2}\mathbb E\left[\sup_{t\in [0,T\wedge\tau_{m}]}(1+\log(1+\|u\|_{H^{s}}^{2}+\|\gamma\|_{H^{s}}^{2}))\bigg|\mathcal{F}_0\right] \nonumber\\
&+C\mathbb E\bigg[\int_{0}^{T\wedge\tau_{m}}\frac{a^{2}(t)(1+\|u\|_{W^{1,\infty}}+\|\gamma\|_{W^{1,\infty}})^{2\theta} (\|T_{\epsilon}u\|_{H^{s}}^{2}+\|T_{\epsilon}\gamma\|_{H^{s}}^{2})^{2}}{(1+\|T_{\epsilon}u\|_{H^{s}}^{2} +\|T_{\epsilon}\gamma\|_{H^{s}}^{2})^{2}(1+\log(1+\|u\|_{H^{s}}^{2}+\|\gamma\|_{H^{s}}^{2}))}d t\bigg|\mathcal{F}_0\bigg]\nonumber\\
&+C_{1}T+C_{2}\mathbb E\bigg[\int_{0}^{T\wedge\tau_{m}}\frac{a^{2}(t)(1+\|u\|_{W^{1,\infty}}+\|\gamma\|_{W^{1,\infty}})^{2\theta}(\|u\|_{H^{s}}^{2}+ \|\gamma\|_{H^{s}}^{2})^{2}}{(1+\|u\|_{H^{s}}^{2}+\|\gamma\|_{H^{s}}^{2})^{2}(1+\log(1+\|u\|_{H^{s}}^{2} +\|\gamma\|_{H^{s}}^{2}))}d t\bigg|\mathcal{F}_0\bigg]\\
&+\mathbb{E}\left[\int_0^{T\wedge\tau_{m}}[|I^\epsilon_1(t)|+|I_2^\epsilon(t)|+|I^\epsilon_3(t)|]d t\bigg|\mathcal{F}_0\right]\nonumber\\
\leq&\frac{1}{2}\mathbb E\left[\sup_{t\in [0,T\wedge\tau_{m}]}(1+\log(1+\|u\|_{H^{s}}^{2}+\|\gamma\|_{H^{s}}^{2}))\bigg|\mathcal{F}_0\right] +\mathbb{E}\left[\int_0^{T\wedge\tau_{m}}[|I^\epsilon_1(t)|+|I_2^\epsilon(t)|+|I^\epsilon_3(t)|]d t\bigg|\mathcal{F}_0\right]\nonumber\\
&+C\mathbb E\bigg[\int_{0}^{T\wedge\tau_{m}}\frac{a^{2}(t)(1+\|u\|_{W^{1,\infty}}+\|\gamma\|_{W^{1,\infty}})^{2\theta} (\|T_{\epsilon}u\|_{H^{s}}^{2}+\|T_{\epsilon}\gamma\|_{H^{s}}^{2})^{2}}{(1+\|T_{\epsilon}u\|_{H^{s}}^{2} +\|T_{\epsilon}\gamma\|_{H^{s}}^{2})^{2}(1+\log(1+\|u\|_{H^{s}}^{2}+\|\gamma\|_{H^{s}}^{2}))}d t\bigg|\mathcal{F}_0\bigg]\\
&+C(u_{0}, \gamma_{0}, C_{1}, C_{2}, T)+C_1T.
\end{align*}
Thus, we use the dominated convergence Theorem, \eqref{lemma3} and \eqref{dom} to obtain
\begin{eqnarray*}
\mathbb E\left[\sup_{t\in [0,T\wedge\tau_{m}]}\log(1+\|u\|_{H^{s}}^{2}+\|\gamma\|_{H^{s}}^{2})\bigg|\mathcal{F}_0\right]\leq C(u_{0}, \gamma_{0}, C_{1}, C_{2}, T).
\end{eqnarray*}
Since $\log(1+x)$ is increasing for $x>0$, we have that for any $m\geq1$,
\begin{align*}
\mathbb P\{\tau_{m}<T|\mathcal{F}_0\}&\leq\mathbb P\bigg\{\sup_{t\in [0,T\wedge \tau_m]}\log(1+\|u\|_{H^{s}}^{2}+\|\gamma\|_{H^{s}}^{2})\geq\log(1+m)\bigg|\mathcal{F}_0\bigg\}\leq\frac{C(u_{0}, \gamma_{0}, C_{1}, C_{2}, T)}{\log(1+m)}.
\end{align*}
Letting $m\rightarrow\infty$ forces $\mathbb P\{\tau^{\ast}<T|\mathcal{F}_0\}=0$ for any $T>0$, which means $\mathbb P\{\tau^{\ast}=\infty\}=1$.

\section{Proof of Theorem \ref{noise I}}
In this section, we study \eqref{SMCH21} with linear noise satisfying Assumption \ref{Linear noise}. Depending on the strength of the noise in \eqref{SMCH21}, we provide the global existence of pathwise solutions for the maximal pathwise solution. Motivated by \cite{T,NV,RZZ}, we introduce
\begin{eqnarray*}
\beta(\omega,t)=e^{\int_{0}^{t}b(t^{'})dW_{t^{'}}-\int_{0}^{t}\frac{b^{2}(t^{'})}{2}dt^{'}}.
\end{eqnarray*}
\begin{pro}\label{6.1}
Let $s>3/2$ and $h_{1}(t, u, \gamma)=b(t)u, h_{2}(t, u, \gamma)=b(t)\gamma$ such that $b(t)$ satisfies Assumption \ref{Linear noise}. Let $(u_{0},\gamma_{0})$ be an $H^{s}\times H^{s}$-valued $\mathcal{F}_{0}$-measurable random variable and $(z,\tau^{\ast})$ be the corresponding unique maximal solution to \eqref{SMCH21}. Let $v_{1}=\beta^{-1}u, v_{2}=\beta^{-1}\gamma$. Then for $t\in[0, \tau^{\ast})$, the processes $v_{1}, v_{2}$ solve the following problem
\begin{eqnarray}\label{v1,v2 2}
\left\{
\begin{array} {l}
\partial_{t}v_{1}+\beta v_{1}v_{1x}+\beta(1-\partial_{x}^{2})^{-1}\partial_{x}(v_{1}^{2}+\frac{1}{2}v_{1x}^{2}+\frac{1}{2}v_{2}^{2}-\frac{1}{2}v_{2x}^{2})=0,\\
\partial_{t}v_{2}+\beta v_{1}v_{2x}+\beta(1-\partial_{x}^{2})^{-1}((v_{1x}v_{2x})_{x}+v_{1x}v_{2})=0,\\
v_{1}(\omega,0,x)=u_{0}(\omega,x), v_{2}(\omega,0,x)=\gamma_{0}(\omega,x).
\end{array}
\right.
\end{eqnarray}
Moreover, we have~$\mathbb P-a.s.$ $(v_{1}, v_{2})\in C([0, \tau^{\ast}); H^{s}\times H^{s})\cap C^{1}([0, \tau^{\ast}); H^{s-1}\times H^{s-1})$. In addition, if $s>5/2$, then it holds
\begin{eqnarray}\label{H1}
\mathbb P\{\|v_{1}\|_{H^{1}}+\|v_{2}\|_{H^{1}}=\|u_{0}\|_{H^{1}}+\|\gamma_{0}\|_{H^{1}}~for~all~t<\tau^\ast\}=1.
\end{eqnarray}
\end{pro}
\begin{pf}
Since $b(t)$ satisfies Assumption \ref{Linear noise}, $h_{1}(t, u, \gamma)=b(t)u, h_{2}(t, u, \gamma)=b(t)\gamma$ satisfy Assumption \ref{h1,h2}, Theorem \ref{Maximal solutions} implies that \eqref{SMCH21} has a unique maximal solution
$(z, \tau^{\ast})$. A direct computation with the It\^o formula yields
\begin{eqnarray*}
d\beta^{-1}=-b(t)\beta^{-1}dW+b^{2}(t)\beta^{-1}dt.
\end{eqnarray*}
Therefore we have
\begin{eqnarray*}
dv_{1}&=&\beta^{-1}du+ud\beta^{-1}+d\beta^{-1}du\nonumber\\
&=&\beta^{-1}\left[-uu_{x}-(1-\partial_{x}^{2})^{-1}\partial_{x}\left(u^{2}+ \frac{1}{2}u_{x}^{2}+\frac{1}{2}\gamma^{2}-\frac{1}{2}\gamma_{x}^{2}\right)\right]dt+b(t)\beta^{-1}udW\nonumber\\
&&+u[-b(t)\beta^{-1}dW+b^{2}(t)\beta^{-1}dt]-b^{2}(t)\beta^{-1}udt\nonumber\\
&=&-\beta v_{1}v_{1x}-\beta(1-\partial_{x}^{2})^{-1}\partial_{x}\left(v_{1}^{2}+\frac{1}{2}v_{1x}^{2}+ \frac{1}{2}v_{2}^{2}-\frac{1}{2}v_{2x}^{2}\right),\nonumber\\
dv_{2}&=&\beta^{-1}d\gamma+\gamma d\beta^{-1}+d\beta^{-1}d\gamma\nonumber\\
&=&\beta^{-1}[-u\gamma_{x}-(1-\partial_{x}^{2})^{-1}((u_{x}\gamma_{x})_{x}+u_{x}\gamma)]dt+b(t)\beta^{-1}\gamma dW\nonumber\\
&&+\gamma[-b(t)\beta^{-1}dW+b^{2}(t)\beta^{-1}dt]-b^{2}(t)\beta^{-1}\gamma dt\nonumber\\
&=&-\beta v_{1}v_{2x}-\beta(1-\partial_{x}^{2})^{-1}((v_{1x}v_{2x})_{x}+v_{1x}v_{2})
\end{eqnarray*}
and since $v_{1}(\omega,0,x)=u_{0}(\omega,x), v_{2}(\omega,0,x)=\gamma_{0}(\omega,x)$, we see that $(v_{1}, v_{2})$ satisfies \eqref{v1,v2 2}. Moreover, Theorem \ref{Maximal solutions}  implies $(u, \gamma) \in C([0, \tau^{\ast}), H^{s}\times H^{s})~\mathbb P-a.s.$, so is $(v_{1}, v_{2})$. From Lemma \ref{F1,F2}
and \eqref{v1,v2 2}, we see that for $\mathbb P-a.s.$, $v_{1t}=-\beta v_{1}v_{1x}-\beta(1-\partial_{x}^{2})^{-1}\partial_{x}(v_{1}^{2}+\frac{1}{2}v_{1x}^{2}+\frac{1}{2}v_{2}^{2}-\frac{1}{2}v_{2x}^{2}), v_{2t}=-\beta v_{1}v_{2x}-\beta(1-\partial_{x}^{2})^{-1}((v_{1x}v_{2x})_{x}+v_{1x}v_{2})$, $(v_{1t}, v_{2t})\in C([0, \tau^{\ast}), H^{s-1}\times H^{s-1})$. Hence, it holds $\mathbb P$-a.s. $(v_{1}, v_{2})\in C^{1}([0, \tau^{\ast}), H^{s-1}\times H^{s-1})$.

In addition, the first two equations of \eqref{v1,v2 2} are equivalent to
\begin{eqnarray}\label{v1,v2 4}
v_{1t}-v_{1xxt}+3\beta v_{1}v_{1x}-2\beta v_{1x}v_{1xx}-\beta v_{1}v_{1xxx}+\beta (v_{2}-v_{2xx})v_{2x}=0,
\end{eqnarray}
\begin{eqnarray}\label{v1,v2 5}
v_{2t}-v_{2xxt}+\beta (v_{1}v_{2x}+v_{1x}v_{2})-\beta (v_{1}v_{2xxx}+v_{1x}v_{2xx})=0.
\end{eqnarray}
Multiplying both sides of \eqref{v1,v2 4} by $v_{1}$ and multiplying both sides of \eqref{v1,v2 5} by $v_{2}$, then integrating the equation on $x\in \mathbb R$, and finally adding the two derived equations, we arrive at
$\mathbb{P}$-a.s.
\begin{eqnarray*}
\frac{d}{dt}\int_{\mathbb R}(v_{1}^{2}+v_{2}^{2}+v_{1x}^{2}+v_{2x}^{2})dx=0,\ \ t<\tau^\ast,
\end{eqnarray*}
which implies \eqref{H1}.
\end{pf}
\textbf{Proof of Theorem 3.7.} To begin with, we apply the operator $D^{s}$ to \eqref{v1,v2 4} and \eqref{v1,v2 5} , multiply both sides of the
resulting equation by $D^{s}v_{1}, D^{s}v_{2}$ respectively, and then integrate on $\mathbb R$ to obtain $\mathbb P-a.s.$
\begin{eqnarray*}
\frac{1}{2}\frac{d}{dt}(\|v_{1}\|_{H^{s}}^{2}+\|v_{2}\|_{H^{s}}^{2})&=&-\beta(\omega,t)\int_{\mathbb R}D^{s}v_{1}D^{s}(v_{1}v_{1x})dx-\beta(\omega,t)\int_{\mathbb R}D^{s}v_{1}D^{s}F_{1}(v_{1},v_{2})dx\nonumber\\
&&-\beta(\omega,t)\int_{\mathbb R}D^{s}v_{2}D^{s}(v_{1}v_{2x})dx-\beta(\omega,t)\int_{\mathbb R}D^{s}v_{2}D^{s}F_{2}(v_{1},v_{2})dx.
\end{eqnarray*}
By \eqref{PPY}, we conclude that $\mathbb P-a.s.$
\begin{eqnarray}\label{*}
\frac{d}{dt}(\|v_{1}\|_{H^{s}}^{2}+\|v_{2}\|_{H^{s}}^{2})\leq C\beta(\omega,t)(\|v_{1}\|_{W^{1,\infty}}+\|v_{2}\|_{W^{1,\infty}})(\|v_{1}\|_{H^{s}}^{2}+\|v_{2}\|_{H^{s}}^{2}).
\end{eqnarray}
Letting $w_{1}=e^{-\int_{0}^{t}b(t^{'})dW_{t^{'}}}u=e^{-\int_{0}^{t}\frac{b^{2}(t^{'})}{2}dt^{'}}v_{1}, w_{2}=e^{-\int_{0}^{t}b(t^{'})dW_{t^{'}}}\gamma=e^{-\int_{0}^{t}\frac{b^{2}(t^{'})}{2}dt^{'}}v_{2}$ and $\alpha(\omega,t)=e^{\int_{0}^{t}b(t^{'})dW_{t^{'}}}$, we obtain
\begin{align*}
&\frac{d}{dt}(\|w_{1}\|_{H^{s}}^{2}+\|w_{2}\|_{H^{s}}^{2})+b^{2}(t)(\|w_{1}\|_{H^{s}}^{2}+\|w_{2}\|_{H^{s}}^{2})\\
&\leq C\alpha(\omega,t)(\|w_{1}\|_{W^{1,\infty}}+\|w_{2}\|_{W^{1,\infty}})(\|w_{1}\|_{H^{s}}^{2}+\|w_{2}\|_{H^{s}}^{2}).
\end{align*}
Assume $\|u_{0}\|_{H^{s}}^{2}+\|\gamma_{0}\|_{H^{s}}^{2}<\frac{b_{\ast}^{2}}{2C^{2}Q^{2}\lambda_{1}^{2}R} <\frac{b_{\ast}^{2}}{C^{2}Q^{2}\lambda_{1}^{2}}$
 and define
\begin{eqnarray}\label{stopping time 1}
\tau_{1}=\inf\bigg\{t<\tau^\ast: \alpha(\omega,t)(\|w_{1}\|_{W^{1,\infty}}+\|w_{2}\|_{W^{1,\infty}})= (\|u\|_{W^{1,\infty}}+\|\gamma\|_{W^{1,\infty}})>\frac{b(t)^2}{C\lambda_{1}}\bigg\}.
\end{eqnarray}
Then it follows from the embedding $(\|u_{0}\|_{W^{1,\infty}}+\|\gamma_{0}\|_{W^{1,\infty}})\leq Q(\|u_{0}\|_{H^{s}}+\|\gamma_{0}\|_{H^{s}})$
that $\mathbb P\{\tau_{1}>0\}=1$, and it holds
\begin{eqnarray*}
\frac{d}{dt}(\|w_{1}\|_{H^{s}}^{2}+\|w_{2}\|_{H^{s}}^{2})+\frac{(\lambda_{1}-1)b^{2}(t)}{\lambda_{1}}(\|w_{1}\|_{H^{s}}^{2}+\|w_{2}\|_{H^{s}}^{2})\leq 0,\ \ t\in [0,\tau_{1}).
\end{eqnarray*}
This implies that for any $0<\lambda_{2}<\frac{\lambda_{1}-1}{\lambda_{1}}$, $\mathbb{P}-a.s.$
\begin{align}\label{6.7}
&\|u(t)\|_{H^{s}}^{2}+\|\gamma(t)\|_{H^{s}}^{2}\nonumber\\
& \leq (\|u_{0}\|_{H^{s}}^{2}+\|\gamma_{0}\|_{H^{s}}^{2})e^{\int_{0}^{t}b(t^{'})dW_{t^{'}}-\int_{0}^{t}\frac{(\lambda_{1}-1)b^{2}(t^{'})}{\lambda_{1}}dt^{'}}\nonumber\\
&=(\|u_{0}\|_{H^{s}}^{2}+\|\gamma_{0}\|_{H^{s}}^{2})e^{\int_{0}^{t}b(t^{'})dW_{t^{'}}-\lambda_{2}\int_{0}^{t}b^{2}(t^{'})dt^{'}}e^{-\int_{0}^{t} \frac{(\lambda_{1}-1)-\lambda_{1}\lambda_{2}}{\lambda_{1}}b^{2}(t^{'})dt^{'}},\ \ t\in [0, \tau_{1}).
\end{align}
Define
\begin{eqnarray*}\label{6.8}
\tau_{2}=\inf\{t>0: e^{\int_{0}^{t}b(t^{'})dW_{t^{'}}-\lambda_{2}\int_{0}^{t}b^{2}(t^{'})dt^{'}}>R \}.
\end{eqnarray*}
Notice that  $\mathbb P\{\tau_{2}>0\}=1$. From \eqref{6.7}, we have
\begin{eqnarray}\label{6.9}
2\|u(t)\|_{H^{s}}^{2}+2\|\gamma(t)\|_{H^{s}}^{2}&\leq& \frac{2b_{\ast}^{2}}{2C^{2}Q^{2}\lambda_{1}^{2}R}\times R\times e^{-\int_{0}^{t}\frac{(\lambda_{1}-1)-\lambda_{1}\lambda_{2}}{\lambda_{1}}b^{2}(t^{'})dt^{'}}\nonumber\\
&=&\frac{b_{\ast}^{2}}{C^{2}Q^{2}\lambda_{1}^{2}}e^{-\int_{0}^{t}\frac{(\lambda_{1}-1)-\lambda_{1}\lambda_{2}}{\lambda_{1}}b^{2}(t^{'})dt^{'}},~~t\in[0, \tau_{1}\wedge\tau_{2}).
\end{eqnarray}
By Assumption \ref{Linear noise}, \eqref{6.9} and \eqref{stopping time 1}, we find that on $[0, \tau_{1}\wedge\tau_{2})$,~$\mathbb P-a.s.$
\begin{eqnarray*}
(\|u\|_{W^{1,\infty}}+\|\gamma\|_{W^{1,\infty}})\leq Q(\|u\|_{H^{s}}+\|\gamma\|_{H^{s}})\leq \frac{b_{\ast}}{C\lambda_{1}}e^{-\int_{0}^{t}\frac{(\lambda_{1}-1)-\lambda_{1}\lambda_{2}}{2\lambda_{1}}b^{2}(t^{'})dt^{'}}\leq \frac{b^{2}(t)}{C\lambda_{1}}e^{-\int_{0}^{t}\frac{(\lambda_{1}-1)-\lambda_{1}\lambda_{2}}{2\lambda_{1}}b^{2}(t^{'})dt^{'}},
\end{eqnarray*}
which together with $\lambda_{2}<\frac{\lambda_{1}-1}{\lambda_{1}}$ and $b^2>0$ derives
\begin{eqnarray*}
\mathbb P\{\tau_{1}\geq\tau_{2}\}=1.
\end{eqnarray*}
Therefore it follows from \eqref{6.9} and Lemma \ref{process control}  that
\begin{eqnarray*}
\mathbb P\{\|u(t)\|_{H^{s}}^{2}+\|\gamma(t)\|_{H^{s}}^{2}\leq \frac{b_{\ast}^{2}}{2C^{2}Q^{2}\lambda_{1}^{2}}~~for~all~t>0\}\geq \mathbb P\{\tau_{2}=\infty\}\geq 1-R^{-2\lambda_{2}},
\end{eqnarray*}
which completes the proof.\\

\section{Proof of Theorem \ref{noise II}}
\subsection{\textbf{Proof of Theorem \ref{noise II}.}} By proposition \ref{6.1}, we can proceed to prove Theorem \ref{noise II}. Since $H^{s}\hookrightarrow C^{2}$ for $s>5/2$, we have $v_{1}, v_{1x}, v_{2}, v_{2x}\in C^{1}([0, \tau^{\ast})\times \mathbb R)$. Then for $x\in\mathbb R$ and $\mathbb P-a.s.$ $\omega\in \Omega$, the problem
\begin{eqnarray}\label{q}
\left\{
\begin{array} {l}
\frac{d q(\omega,t,x)}{dt}=\beta(\omega,t)v_{1}(\omega,t,q(\omega,t,x)), ~~t\in[0, \tau^{\ast}),\\
q(\omega,0,x)=x,~~x\in\mathbb R
\end{array}
\right.
\end{eqnarray}
has a unique solution $q(\omega, t, x)$ such that $q(\omega, t, x)\in C^{1}([0, \tau^{\ast})\times\mathbb R)$ for $\mathbb P-a.s.$ $\omega\in \Omega$. Moreover, differentiating
\eqref{q} with respect to $x$ yields that for $\mathbb P-a.s.$ $\omega\in\Omega$,
\begin{eqnarray}\label{q1}
\left\{
\begin{array} {l}
\frac{d q_{x}(\omega,t,x)}{dt}=\beta(\omega,t)v_{1x}(\omega,t,q(\omega,t,x))q_{x}, ~~t\in[0, \tau^{\ast}),\nonumber\\
q_x(\omega,0,x)=1,~~x\in\mathbb R.
\end{array}
\right.
\end{eqnarray}
For $\mathbb P-a.s.$ $\omega\in\Omega$, we solve the above equation to obtain
\begin{eqnarray*}
q_{x}(\omega,t,x)=\exp\bigg(\int_{0}^{t}\beta(\omega,t^{'})v_{1x}(\omega,t^{'},q(\omega,t^{'},x))dt^{'}\bigg), \ \ t\in(0,\tau^\ast).
\end{eqnarray*}
Thus for $\mathbb P-a.s.$ $\omega\in\Omega$, $q_{x}(\omega,t,x)>0$, $(t, x)\in [0, \tau^{\ast})\times\mathbb R$. Then the momentum variable $V_{1}=(1-\partial_{x}^{2})v_{1}, V_{2}=(1-\partial_{x}^{2})v_{2}$ satisfy $\mathbb P-a.s.$
\begin{eqnarray}\label{v11}
V_{1t}+\beta v_{1}V_{1x}+2\beta v_{1x}V_{1}+\beta  v_{2x}V_{2}=0,\nonumber\\
V_{2t}++\beta (v_{1}V_{2})_{x}=0.
\end{eqnarray}
Applying particle trajectory method \eqref{q} and the first equation of \eqref{v11}, we obtain
\begin{align*}
&\frac{d}{dt}\left[e^{\int_{0}^{t}\frac{\beta (\omega,s)V_{2}v_{2x}(\omega,s,x)}{V_{1}(\omega,s,x)}ds}V_{1}(\omega,t,q(\omega,t,x))q_{x}^{2}(\omega,t,x)\right]\nonumber\\
=&e^{\int_{0}^{t}\frac{\beta (\omega,s)V_{2}v_{2x}(\omega,s,x)}{V_{1}(\omega,s,x)}ds} \beta (\omega,s)q_{x}^{2}V_{2}v_{2x}(\omega,s,x)
+e^{\int_{0}^{t}\frac{\beta (\omega,s)V_{2}v_{2x}(\omega,s,x)}{V_{1}(\omega,s,x)}ds}q_{x}^{2}[V_{1t}+\beta v_{1}V_{1x}+2\beta v_{1x}V_{1}]\nonumber\\
=&e^{\int_{0}^{t}\frac{\beta (\omega,s)V_{2}v_{2x}(\omega,s,x)}{V_{1}(\omega,s,x)}ds}q_{x}^{2}[V_{1t}+\beta v_{1}V_{1x}+2\beta v_{1x}V_{1}+\beta  v_{2x}V_{2}]\nonumber\\
=&0.
\end{align*}
This and $q_{x}(\omega,0,x)=1$ imply that
\begin{eqnarray}\label{v1 q 1}
e^{\int_{0}^{t}\frac{\beta (\omega,s)V_{2}v_{2x}(\omega,s,x)}{V_{1}(\omega,s,x)}ds}V_{1}(\omega,t,q(\omega,t,x))q_{x}^{2}(\omega,t,x)=V_{1}(\omega,0,x).
\end{eqnarray}
Consequently, we have sign$(V_{1}(\omega,t,x))$=sign$(V_{1}(\omega,0,x))$.

The next step, we give the following useful lemma that will be used in the sequel.
\begin{lem}\label{blow t}(Blow-up scenario 2)
Let  $s>3/2$ and $(u_{0}, \gamma_{0})$ be an $H^{s}\times H^{s}$-valued $\mathcal{F}_{0}$-measurable random variable. Assume that $(z,\tau^{\ast})$ is the corresponding maximal solution. Then $z$ as a $W^{1,\infty}\times W^{1,\infty}$-valued process is $\mathcal{F}_{t}$-adapted for $t<\tau^{\ast}$ and $\mathbb P-a.s.$ on the set $\{\tau^\ast<\infty\}$
\begin{eqnarray}\label{blow}
&&1_{\{\lim\mathop{\sup}\limits_{t\rightarrow \tau^{\ast}}(\|u(t)\|_{H^{s}}+\|\gamma(t)\|_{H^{s}})=\infty\}}=1_{\{\lim\mathop{\sup}\limits_{t\rightarrow \tau^{\ast}}\|u(t)\|_{W^{1,\infty}}=\infty\}}.
\end{eqnarray}
\begin{pf}
It is clear that $\{\lim\mathop{\sup}\limits_{t\rightarrow \tau^{\ast}}\|u(t)\|_{W^{1,\infty}}=\infty\}\subset\{\lim\mathop{\sup}\limits_{t\rightarrow \tau^{\ast}}(\|u(t)\|_{H^{s}}+\|\gamma(t)\|_{H^{s}})=\infty\}$. It is sufficient to prove  $\{\lim\mathop{\sup}\limits_{t\rightarrow \tau^{\ast}}\|u(t)\|_{W^{1,\infty}}=\infty\}^{C}\subset\{\lim\mathop{\sup}\limits_{t\rightarrow \tau^{\ast}}(\|u(t)\|_{H^{s}}+\|\gamma(t)\|_{H^{s}})=\infty\}^{C}$.
Notice that
\begin{eqnarray}\label{***}
\{\lim\mathop{\sup}\limits_{t\rightarrow \tau^{\ast}}\|u(\omega,t)\|_{W^{1,\infty}}=\infty\}^{C}=\{\exists M(\omega)>0, s.t.~ \|u(\omega,t)\|_{W^{1,\infty}}\leq M(\omega),~~\forall t <\tau^\ast\}.
\end{eqnarray}
By the equation \eqref{v1,v2 5} and using the identity $\partial_{x}^{2}G\ast f=\partial_{x}^{2}(1-\partial_{x}^{2})^{-1}f=(1-\partial_{x}^{2})^{-1}f-f$, we have
\begin{align}\label{TTY}
&\bigg|\frac{dv_{2x}(\omega,t,q(\omega,t,x))}{dt}\bigg|\nonumber\\
&=|v_{2tx}(t,q)+v_{2xx}(t,q)\beta v_{1}|\nonumber\\
&=|-\beta v_{1x}v_{2x}-\beta\partial_{x}^{2}(1-\partial_{x}^{2})^{-1}(v_{1x}v_{2x})-\beta\partial_{x}(1-\partial_{x}^{2})^{-1}(v_{1x}v_{2})|\nonumber\\
&=|\beta(1-\partial_{x}^{2})^{-1}(v_{1x}v_{2x})-\beta\partial_{x}(1-\partial_{x}^{2})^{-1}(v_{1x}v_{2})|\nonumber\\
&\leq\beta\|G\|_{L^{\infty}}\|v_{1x}v_{2x}\|_{L^{1}}+\beta\|\partial_{x}G\|_{L^{\infty}}\|v_{1x}v_{2}\|_{L^{1}}\nonumber\\
&\leq C\beta(2\|v_{1x}\|_{L^{2}}+\|v_{2x}\|_{L^{2}}+\|v_{2}\|_{L^{2}})\leq C\beta(\|v_{1}(0)\|_{H^{1}}+\|v_{2}(0)\|_{H^{1}}).
\end{align}
For $m\geq 1$, define  $$\tau_m=\inf\{t<\tau^\ast:\|u(t)\|_{H^{s}}+\|\gamma(t)\|_{H^{s}}\geq m\}.$$
By \eqref{TTY}, Sobolev's embedding and \eqref{H1}, we have
\begin{eqnarray}\label{gamma}
\|v_{2}(\omega,t,q(\omega,t,\cdot))\|_{W^{1,\infty}}\leq C \int_{0}^{t}\beta(\omega,t^{'}) dt^{'} (\|u_{0}\|_{H^{1}}+\|\gamma_{0}\|_{H^{1}})+\|\gamma_{0}\|_{W^{1,\infty}},\ \ t\leq \tau_m.
\end{eqnarray}
In addition, we derive from \eqref{*} that
\begin{eqnarray*}
\frac{d}{dt}(\|v_{1}\|_{H^{s}}^{2}+\|v_{2}\|_{H^{s}}^{2})\leq C(\|u\|_{W^{1,\infty}}+\beta(\omega,t)\|v_{2}\|_{W^{1,\infty}})(\|v_{1}\|_{H^{s}}^{2}+\|v_{2}\|_{H^{s}}^{2}).
\end{eqnarray*}
By means of Gronwall's inequality and \eqref{***}, for any $\omega\in\{\lim\mathop{\sup}\limits_{t\rightarrow \tau^{\ast}}\|u(\omega,t)\|_{W^{1,\infty}}=\infty\}^{C}$, we obtain
\begin{align*}
&\|v_{1}(T\wedge \tau_m)\|_{H^{s}}^{2}+\|v_{2}(T\wedge \tau_m)\|_{H^{s}}^{2} \nonumber\\
&\leq (\|u_{0}\|_{H^{s}}^{2}+\|\gamma_{0}\|_{H^{s}}^{2}) \exp\left\{\int_0^{T\wedge\tau_m}C(\|u\|_{W^{1,\infty}}+\beta(\omega,t)\|v_{2}\|_{W^{1,\infty}})dt\right\},\nonumber\\
&\leq (\|u_{0}\|_{H^{s}}^{2}+\|\gamma_{0}\|_{H^{s}}^{2})\\
&\times \exp\left\{C \bigg(M(T\wedge \tau_m)+\int_{0}^{T\wedge \tau_m}\beta(\omega,t)\bigg[\int_{0}^{t}\beta(\omega,t^{'}) dt^{'} (\|u_{0}\|_{H^{1}}+\|\gamma_{0}\|_{H^{1}})+\|\gamma_{0}\|_{W^{1,\infty}}\bigg]dt\bigg)\right\}.\nonumber
\end{align*}
This implies on the set $\{\tau^\ast<\infty\}\cap\{\lim\mathop{\sup}\limits_{t\rightarrow \tau^{\ast}}\|u(\omega,t)\|_{W^{1,\infty}}=\infty\}^{C}$,
\begin{align*}
&\|u(T\wedge \tau_m)\|_{H^{s}}^{2}+\|\gamma(T\wedge \tau_m)\|_{H^{s}}^{2}\\
\leq& (\|u_{0}\|_{H^{s}}^{2}+\|\gamma_{0}\|_{H^{s}}^{2})\beta(\omega,T\wedge \tau_m)\\
&\times  \exp\left\{C \bigg(M\tau_m+\int_{0}^{T\wedge\tau_m}\bigg[\beta(\omega,t)\int_{0}^{t}\beta(\omega,t^{'}) dt^{'} (\|u_{0}\|_{H^{1}}+\|\gamma_{0}\|_{H^{1}})+\|\gamma_{0}\|_{W^{1,\infty}}\bigg]dt\bigg)\right\}\\
<&\infty,
\end{align*}
where we used $\sup_{t>0}\beta(\omega,t)<\infty$ due to $\sup_{t>0}\mathbb E \beta(\omega, t)=1$ and Doob's $L^{1}$-inequality. Hence we can see that on the set $\{\tau^\ast<\infty\}$, $\{\lim\mathop{\sup}\limits_{t\rightarrow \tau^{\ast}}\|u(t)\|_{W^{1,\infty}}=\infty\}^{C}\subset\{\lim\mathop{\sup}\limits_{t\rightarrow \tau^{\ast}}(\|u(t)\|_{H^{s}}+\|\gamma(t)\|_{H^{s}})=\infty\}^{C}$. So, we finish the proof.
\end{pf}
\end{lem}

\begin{lem}\label{blow t2}(Blow-up scenario 3)
Let $s>3/2$ and $z_{0}$ be an $H^{s}\times H^{s}$-valued $\mathcal{F}_{0}$-measurable random variable. Assume that $(z,\tau^{\ast})$ is the corresponding maximal solution. Then $z$ as a $W^{1,\infty}\times W^{1,\infty}$-valued process is $\mathcal{F}_{t}$-adapted for $t<\tau^{\ast}$ and $\mathbb P-a.s.$ on the set $\{\tau^\ast<\infty\}$,
\begin{eqnarray}\label{blow 1}
&&1_{\{\lim\mathop{\sup}\limits_{t\rightarrow \tau^{\ast}}(\|u(t)\|_{H^{s}}+\|\gamma(t)\|_{H^{s}})=\infty\}}=1_{\{\lim\inf_{t\rightarrow \tau^{\ast}}\min_{x\in\R}\{u_{x}(\omega, t, x)\}=-\infty\}}.
\end{eqnarray}
\end{lem}
\begin{pf} It is clear that $\{\lim\inf_{t\rightarrow \tau^{\ast}}\min_{x\in\R}\{u_{x}(\omega, t, x)\}=-\infty\}\subset\{\lim\mathop{\sup}\limits_{t\rightarrow \tau^{\ast}}(\|u(t)\|_{H^{s}}+\|\gamma(t)\|_{H^{s}})=\infty\}$.
The rest of proof is similar to that of Lemma \ref{blow t} by replacing equation \eqref{***} with
\begin{align}\label{TTW}
\{\lim\inf_{t\rightarrow \tau^{\ast}}\min_{x\in\R}\{u_{x}(\omega, t, x)\}=-\infty\}^{C}=\{\exists M(\omega)>0, s.t.~u_{x}(\omega, t, x)>-M(\omega),~~\forall t<\tau^\ast\}.
\end{align}
Without loss of generality, we only need to show that this Lemma holds for $s=2$. Multiplying the first equation in \eqref{v11} by $V_{1}=(1-\partial_{x}^{2})v_{1}$ and integrating by parts, we get
\begin{eqnarray}\label{V1}
\frac{d}{dt}\int_{\mathbb R}V_{1}^{2}dx&=&-2\beta(w,t)\int_{\mathbb R}v_{1}V_{1}V_{1x}dx-4\beta(w,t)\int_{\mathbb R}V_{1}^{2}v_{1x}dx-2\beta(w,t)\int_{\mathbb R}V_{1}V_{2}v_{2x}dx\nonumber\\
&=&-3\beta(w,t)\int_{\mathbb R}V_{1}^{2}v_{1x}dx-2\beta(w,t)\int_{\mathbb R}V_{1}V_{2}v_{2x}dx.
\end{eqnarray}
Multiplying the second equation in \eqref{v11} by $V_{2}=(1-\partial_{x}^{2})v_{2}$ and integrating by parts, we obtain
\begin{eqnarray}\label{V2}
\frac{d}{dt}\int_{\mathbb R}V_{2}^{2}dx&=&-\beta(w,t)\int_{\mathbb R}v_{1x}V_{2}^{2}dx.
\end{eqnarray}
Thus, in view of \eqref{TTW}, \eqref{V1}, \eqref{V2} and \eqref{gamma}, for any $\omega\in \{\lim\inf_{t\rightarrow \tau^{\ast}}\min_{x\in\R}\{u_{x}(\omega, t, x)\}=-\infty\}^{C}$, we obtain
\begin{align*}
&\frac{d}{dt}\int_{\mathbb R}(V_{1}^{2}+V_{2}^{2})dx\nonumber\\
&=-3\beta(w,t)\int_{\mathbb R}V_{1}^{2}v_{1x}dx-\beta(w,t)\int_{\mathbb R}v_{1x}V_{2}^{2}dx-2\beta(w,t)\int_{\mathbb R}V_{1}V_{2}v_{2x}dx\nonumber\\
&\leq 3M\int_{\mathbb R}(V_{1}^{2}+V_{2}^{2})dx\nonumber\\
&+\beta(\omega,t)\bigg(\int_{0}^{t}\beta(\omega,t^{'}) dt^{'} (\|u_{0}\|_{H^{1}}+\|\gamma_{0}\|_{H^{1}})+\|\gamma_{0}\|_{W^{1,\infty}}\bigg)\int_{\mathbb R}(V_{1}^{2}+V_{2}^{2})dx
\end{align*}
By means of Gronwall's inequality, we arrive at
\begin{eqnarray}\label{V1V20}
&&\|v_{1}(T\wedge\tau_m)\|_{H^{2}}+\|v_{2}(T\wedge\tau_m)\|_{H^{2}}=\|V_{1}(T\wedge\tau_m)\|_{L^{2}}+ \|V_{2}(T\wedge\tau_m)\|_{L^{2}}\nonumber\\
&&\leq(\|V_{1}(0,\cdot)\|_{L^{2}}+\|V_{2}(0,\cdot)\|_{L^{2}})\nonumber\\
&&\times \exp\bigg[\int_{0}^{T\wedge\tau_m}\left(3M+\beta(w,t)\bigg(\int_{0}^{t}\beta(\omega,t^{'}) dt^{'} (\|u_{0}\|_{H^{1}}+\|\gamma_{0}\|_{H^{1}})+\|\gamma_{0}\|_{W^{1,\infty}}\bigg)\right)dt\bigg].\nonumber
\end{eqnarray}
 Then on the set $\{\tau^\ast<\infty\}\cap\{\lim\inf_{t\rightarrow \tau^{\ast}}\min_{x\in\R}\{u_{x}(\omega, t, x)\}=-\infty\}^{C}$,
\begin{align*}
&\|u(T\wedge\tau_m)\|_{H^{2}}+\|\gamma(T\wedge\tau_m)\|_{H^{2}}\nonumber\\
&\leq\beta(w,T\wedge\tau_m)(\|V_{1}(0,\cdot)\|_{L^{2}}+\|V_{2}(0,\cdot)\|_{L^{2}})\nonumber\\
&\times\exp\bigg[\int_{0}^{T\wedge\tau_m}\left(3M+\beta(w,t)\bigg(\int_{0}^{t}\beta(\omega,t^{'}) dt^{'} (\|u_{0}\|_{H^{1}}+\|\gamma_{0}\|_{H^{1}})+\|\gamma_{0}\|_{W^{1,\infty}}\bigg)\right)dt\bigg]\nonumber\\
&<\infty,
\end{align*}
where we used $\sup_{t>0}\beta(\omega,t)<\infty$ due to $\sup_{t>0}\mathbb E \beta(\omega, t)=1$ and Doob's $L^{1}$-inequality. Hence we can see that $\{\lim\inf_{t\rightarrow \tau^{\ast}}\min_{x\in\R}\{u_{x}(\omega, t, x)\}=-\infty\}^{C}\subset\{\lim\mathop{\sup}\limits_{t\rightarrow \tau^{\ast}}(\|u(t)\|_{H^{s}}+\|\gamma(t)\|_{H^{s}})=\infty\}^{C}$.
This completes the proof.
\end{pf}
\begin{lem}\label{lem 6.3}
Let $V_{1}=v_{1}-v_{1xx}$ and $\mathbb P\{V_{1}(\omega, 0, x)>0, \forall x\in \mathbb R\}=p$, $\mathbb P\{V_{1}(\omega, 0, x)<0, \forall x\in \mathbb R\}=q$, $\mathbb P\bigg\{V_{1}(\omega, 0, x)\leq0, ~x\leq x_{0} ~~and~~ V_{1}(\omega, 0, x)\geq0, ~x\geq x_{0}\bigg\}=m$ for some $p, q, m\in[0, 1]$. Then the maximal solution $(z,\tau^{\ast})$ of \eqref{SMCH21} satisfies
\begin{eqnarray}\label{P}
&&\mathbb P\bigg\{\|u_{x}(\omega, t)\|_{L^{\infty}}\leq \frac{\sqrt{2}}{2}\beta(\omega, t)(\|u_{0}\|_{H^{1}}+\|\gamma_{0}\|_{H^{1}})~~ \nonumber\\
&&~or~~ u_{x}(\omega, t)\geq-\frac{\sqrt{2}}{2}\beta(\omega, t,x)(\|u_{0}\|_{H^{1}}+\|\gamma_{0}\|_{H^{1}}),~~\forall t\in [0,\tau^{\ast})\bigg\}\geq p+q+m.
\end{eqnarray}
\end{lem}

\begin{pf} Denote
\begin{align*}
&A_p=\{V_{1}(\omega, 0, x)>0, \forall x\in \mathbb R\},\\
&A_q=\{V_{1}(\omega, 0, x)<0, \forall x\in \mathbb R\},\\
&A_m=\{V_{1}(\omega, 0, x)\leq0, ~x\leq x_{0} ~~and~~ V_{1}(\omega, 0, x)\geq0, ~x\geq x_{0}\}.
\end{align*}
Using $G(x)=\frac{e^{-|x|}}{2}$, one can derive that for $\mathbb{P}$-a.s. $\omega\in\Omega$, and for all $(t,x)\in [0, \tau^{\ast})\times \mathbb R$,
\begin{eqnarray*}
v_{1}(\omega,t,x)=\frac{1}{2}e^{-x}\int_{-\infty}^{x}e^{\xi}V_{1}(\omega,t,\xi)d\xi+\frac{1}{2}e^{x}\int_{x}^{\infty}e^{-\xi}V_{1}(\omega,t,\xi)d\xi,
\end{eqnarray*}
\begin{eqnarray*}
v_{1x}(\omega,t,x)=-\frac{1}{2}e^{-x}\int_{-\infty}^{x}e^{\xi}V_{1}(\omega,t,\xi)d\xi+\frac{1}{2}e^{x}\int_{x}^{\infty}e^{-\xi}V_{1}(\omega,t,\xi)d\xi.
\end{eqnarray*}
Therefore,
\begin{eqnarray}\label{u ux 3}
v_{1}(\omega,t,x)+v_{1x}(\omega,t,x)=e^{x}\int_{x}^{\infty}e^{-\xi}V_{1}(\omega,t,\xi)d\xi,
\end{eqnarray}
\begin{eqnarray}\label{u ux 4}
v_{1}(\omega,t,x)-v_{1x}(\omega,t,x)=e^{-x}\int_{-\infty}^{x}e^{\xi}V_{1}(\omega,t,\xi)d\xi.
\end{eqnarray}
Then one can employ \eqref{u ux 3}, \eqref{u ux 4} and sign$(V_{1})$=sign$(V_{1}(\omega,0,x))$, $(t,x)\in [0, \tau_{\ast})\times \mathbb R$ to obtain that for all $(t, x)\in [0, \tau^{\ast})\times \mathbb R$,
\begin{eqnarray}\label{v vx}
\left\{
\begin{array} {l}
-v_{1}(\omega,t,x)\leq v_{1x}(\omega,t,x)\leq v_{1}(\omega,t,x),\ \ \omega\in A_p,\\
v_{1}(\omega,t,x)\leq v_{1x}(\omega,t,x)\leq -v_{1}(\omega,t,x),\ \ \omega\in A_q.
\end{array}
\right.
\end{eqnarray}
In addition, since $q(\omega, t, \cdot)$ is an increasing diffeomorphism of $\mathbb R $ with $q_{x}(\omega, t, x)>0$ for all $(t, x)\in[0, \tau^{\ast})\times\mathbb R$, by \eqref{v1 q 1}, it follows that for any $\omega\in A_m$,
\begin{eqnarray}\label{V1 q}
\left\{
\begin{array} {l}
V_{1}(\omega,t,x)\leq 0~~if~~ x\leq q(\omega,t,x_{0}),\\
V_{1}(\omega,t,x)\geq 0~~if~~ x\geq q(\omega,t,x_{0}),\\
V_{1}(\omega,t,q(\omega,t,x_{0}))=0.
\end{array}
\right.
\end{eqnarray}
Therefore, for any $\omega\in A_m$, when $x\leq q(\omega,t,x_{0})$, by \eqref{u ux 4} and \eqref{V1 q}, we have $v_{1}(\omega,t,x)\leq v_{1x}(\omega,t,x)$; when $x\geq q(\omega,t,x_{0})$, by \eqref{u ux 3} and \eqref{V1 q}, we have $v_{1x}(\omega,t,x)\geq-v_{1}(\omega,t,x)$, Therefore, it follows from \eqref{H1} that for any $\omega\in A_m$,
\begin{eqnarray}\label{V01}
-v_{1x}(\omega, t, x)\leq|v_{1}(\omega, t, x)|\leq\|v_{1}(\omega, t, x)\|_{L^{\infty}}\leq\frac{\sqrt{2}}{2}(\|u_{0}\|_{H^{1}}+\|\gamma_{0}\|_{H^{1}}), ~\forall (t,x)\in [0,\tau^{\ast})\times\mathbb R
\end{eqnarray}
Then for any $\omega\in A_m$, $u_{x}(\omega, t)\geq-\frac{\sqrt{2}}{2}\beta(\omega, t,x)(\|u_{0}\|_{H^{1}}+\|\gamma_{0}\|_{H^{1}})$. This together with Lemma 7.2 and $\sup_{t>0}\beta(\omega,t,x)<\infty$ implies that $z$ globally exists.\\

For any $\omega \in A_p\cup A_q $, it follows from \eqref{v vx} that $|v_{1x}(\omega, t, x)|\leq|v_{1}(\omega, t, x)|$, in view of Sobolev inequality and \eqref{H1}, we arrive at
\begin{align}\label{P2}
\|v_{1x}(\omega, t, x)\|_{L^{\infty}}\leq\|v_{1}(\omega, t, x)\|_{L^{\infty}}
\leq&\frac{\sqrt{2}}{2}(\|u_{0}\|_{H^{1}}+\|\gamma_{0}\|_{H^{1}}), ~\forall (t,x)\in [0,\tau^{\ast})\times\mathbb R,
\omega \in A_p\cup A_q .\end{align}

Combining $A_p\cap A_q \cap A_m=\emptyset$, \eqref{V01} and \eqref{P2}, we derive \eqref{P}.
%
\end{pf}
\textbf{Proof of Theorem \ref{noise II}.} Note that $\sup_{t>0}\mathbb E\beta(\omega,t,x)=1$ and Doob's $L^{1}$-inequality implies that $\sup_{t>0}\beta(\omega,t,x)<\infty$. Then we can infer from \eqref{blow},  \eqref{blow 1} and \eqref{P} that $\mathbb P\{\tau^{\ast}=\infty\}\geq p+q+m$. This completes the proof.
\subsection{Proof of Theorem \ref{blow up initial value}}
The proof of Theorem  \ref{blow up initial value} relies on certain properties of the solution $v_{1}, v_{2}$ to the equations \eqref{v1,v2 4} and \eqref{v1,v2 5}. We first prove the following lemma.
\begin{lem}\label{p6.2}
Let $s>5/2$ and $b(t)$ satisfy Assumption \ref{Linear noise}. Assume $(u_{0}, \gamma_{0})$ is an $H^{s}\times H^{s}$-valued $\mathcal{F}_{0}$-measurable random variable. Let $K=\frac{\sqrt{2}}{2} (\|u_{0}\|_{H^{1}}^{2}+\|\gamma_{0}\|_{H^{1}}^{2})^{\frac{1}{2}}$. Then for $v_{1}, v_{2}$ defined by \eqref{v1,v2 4}, \eqref{v1,v2 5} and any $x_0\in\R$,
\begin{eqnarray*}
g(\omega,t):=v_{1x}(\omega,t, q(\omega, t,x_{0}))
\end{eqnarray*}
satisfies $\mathbb P-a.s.$
\begin{eqnarray}\label{6.19}
\frac{d}{dt}g(t)\leq \beta K^{2}-\frac{\beta}{2}g^{2}(t),
 \ \ t<\tau^\ast.
\end{eqnarray}
Moreover, if there exists some $x_{0}\in \mathbb R$ such that $\mathbb P-a.s.$ $g(0)<-\sqrt{2}K$, then $\mathbb P-a.s.$. $g(t)$ is non-increasing on $[0, \tau^{\ast})$ and
\begin{eqnarray}\label{6.20}
g(t)<-\sqrt{2}K, \ \ t\in[0, \tau^{\ast}).
\end{eqnarray}
\end{lem}
\begin{pf}
For any $v_{1}, v_{2}\in H^{1}$, by the  representation of $G*f=(1-\partial_{x}^{2})^{-1}f$, we have
\begin{eqnarray}\label{G1}
G*\bigg(v_{1}^{2}+\frac{1}{2}v_{1x}^{2}\bigg)(x)=\frac{1}{2}\int_{-\infty}^{x}e^{-x+y}\bigg(v_{1}^{2}+\frac{1}{2}v_{1x}^{2}\bigg)(y)dy+\frac{1}{2}\int_{x}^{\infty}e^{x-y}\bigg(v_{1}^{2}+\frac{1}{2}v_{1x}^{2}\bigg)(y)dy.
\end{eqnarray}
The following inequality
\begin{eqnarray*}
\int_{-\infty}^{x}e^{y}\bigg(v_{1}^{2}+v_{1x}^{2}\bigg)(y)dy\geq2\int_{-\infty}^{x}e^{y} v_{1}v_{1x}(y)dy=e^{x}v_{1}^{2}(x)-\int_{-\infty}^{x}e^{y}v_{1}^{2}dy
\end{eqnarray*}
implies that
\begin{eqnarray}\label{G2}
\frac{1}{2}\int_{-\infty}^{x}e^{-x+y}\bigg(v_{1}^{2}+\frac{1}{2}v_{1x}^{2}\bigg)(y)dy\geq \frac{1}{4}v_{1}^{2}(x).
\end{eqnarray}
Similarly, we get the estimate of the second term in \eqref{G1} as
\begin{eqnarray}\label{G3}
\frac{1}{2}\int_{x}^{\infty}e^{x-y}\bigg(v_{1}^{2}+\frac{1}{2}v_{1x}^{2}\bigg)(y)dy\geq \frac{1}{4}v_{1}^{2}(x),
\end{eqnarray}
Combining \eqref{G1}, \eqref{G2} and \eqref{G3}, we deduce $G*(v_{1}^{2}+\frac{1}{2}v_{1x}^{2})(x)\geq \frac{1}{2}v_{1}^{2}(x)$. In addition,
\begin{eqnarray}\label{G4}
\|G*v_{2x}^{2}\|_{L^{\infty}}\leq\|G\|_{L^{\infty}}\|v_{2x}^{2}\|_{L^{1}}=\frac{1}{2}\|v_{2x}^{2}\|_{L^{1}}.
\end{eqnarray}
Differentiating the first equation of \eqref{v1,v2 2} with respect to $x$, and using \eqref{H1} and \eqref{G4}, we have
\begin{eqnarray*}
\frac{d}{dt}v_{1x}(\omega,t,q(t,\omega,x))&=&v_{1xt}+v_{1xx}\beta(\omega,t,x)v_{1}(\omega,t,q(\omega,t,x))\nonumber\\
&=&-\beta v_{1x}^{2}-\beta\partial_{x}^{2}(1-\partial_{x}^{2})^{-1}\left(v_{1}^{2}+\frac{1}{2}v_{1x}^{2}+ \frac{1}{2}v_{2}^{2}-\frac{1}{2}v_{2x}^{2}\right)\nonumber\\
&\leq&-\frac{1}{2}\beta v_{1x}^{2}+\frac{1}{2}\beta v_{1}^{2}+\frac{1}{4}\beta v_{2}^{2}+\frac{3}{4}\beta G*(v_{2x}^{2})\nonumber\\
&\leq&-\frac{1}{2}\beta v_{1x}^{2}+\frac{\beta}{2}(\|u_{0}\|_{H^{1}}^{2}+\|\gamma_{0}\|_{H^{1}}^{2}).
\end{eqnarray*}
In view of the assumptions of Lemma \ref{p6.2}, we have~$\mathbb P-a.s.$
\begin{eqnarray*}
\frac{d}{dt}g(t)\leq-\frac{\beta}{2}g^{2}(t)+\beta K^{2},\ \ t<\tau^\ast,
\end{eqnarray*}
which is \eqref{6.19}. In order to prove \eqref{6.20}, define
\begin{eqnarray*}
\zeta(w):=\inf\bigg\{t\in[0,\tau^\ast): g(w, t)>-\sqrt{2}K\bigg\}.
\end{eqnarray*}
If $g(0)<-\sqrt{2}K$, then $\mathbb P\{\zeta>0\}=1$. From the definition of $\zeta(\omega)$, we find that $\zeta(\omega)\leq\tau^{*}, $ for $\mathbb P-a.s.$ $w\in\Omega$. From \eqref{6.19}, we have that $g(\omega, t)$ is nonincreasing for $t\in[0, \zeta(\omega))$. Hence by the continuity of the path of $g(\omega, t)$, we obtain that $g(\omega, t)\leq g(0)<-\sqrt{2}K,\ \ t\in[0,\zeta(\omega))$. In view of the time continuity of $g(\omega, t)$ again, we find that $\mathbb P\{\zeta=\tau^{*}\}=1.$ Hence \eqref{6.20} is true.

\end{pf}
\textbf{Proof of Theorem \ref{blow up initial value}.} From Lemma \ref{p6.2} and \eqref{u0x}, we rewrite \eqref{6.19} as
\begin{eqnarray*}
\frac{d}{dt}g(t)&\leq& -\frac{\beta(t)}{2}\bigg(1-\frac{2K^{2}}{g^{2}(0)}\bigg)g^{2}(t)-\left(\frac{g^{2}(t)}{g^{2}(0)}-1\right)\beta(t)K^{2}\nonumber\\
&\leq& -\frac{\beta(t)}{2}\bigg(1-\frac{2K^{2}}{g^{2}(0)}\bigg)g^{2}(t),\ \ t\in[0, \tau^{\ast}).
\end{eqnarray*}
Integrating on both sides leads to $\mathbb P-a.s.$
\begin{eqnarray*}
\frac{1}{g(t)}-\frac{1}{g(0)}\geq \left(1-\frac{2K^{2}}{g^{2}(0)}\right)\int_{0}^{t}\frac{\beta(t^{'})}{2}dt^{'},~~t\in[0, \tau^{*}).
\end{eqnarray*}
Assuming $\Omega^{'}=\{\omega: \beta(t,\omega)\geq c e^{-\frac{b^{*}}{2}t}~for~ all~ t\}$, $g(t)\leq-\sqrt{2}K$ means that $\mathbb{P}-a.s.$ $\omega\in\Omega^{'}$
\begin{eqnarray*}
-\frac{1}{g(0)}\geq \bigg(\frac{1}{2}-\frac{K^{2}}{g^{2}(0)}\bigg)\int_{0}^{\tau^{*}}\beta(t^{'})dt^{'}\geq\bigg(\frac{1}{2}-\frac{K^{2}}{g^{2}(0)}\bigg)\bigg(\frac{2c}{b^{*}}-\frac{2c}{b^{*}}e^{-\frac{b^{*}}{2}\tau^{*}}\bigg).
\end{eqnarray*}
If  $g(0)<-\frac{1}{2}\sqrt{\frac{(b^{*})^{2}}{c^{2}}+8 K^{2}}-\frac{b^{*}}{2c}$, we obtain on $\Omega'$
\begin{eqnarray*}
\bigg(\frac{1}{2}-\frac{K^{2}}{g^{2}(0)}\bigg)\frac{2c}{b^{*}}e^{-\frac{b^{*}}{2}\tau^{*}} \geq\frac{2c}{b^{*}}\bigg(\frac{1}{2}-\frac{K^{2}}{g^{2}(0)}\bigg)+\frac{1}{g(0)}>0.
\end{eqnarray*}
Therefore we have $\tau^{*}<\infty$ $\mathbb{P}-a.s.$ on $\Omega^{'}$, which implies that
\begin{eqnarray*}
\mathbb P\{\tau^{*}<\infty\}\geq\mathbb P\{\beta(t)\geq c e^{-\frac{b^{*}}{2}t}~for~ all~ t\}=\mathbb P\left\{e^{\int_{0}^{t}b(t^{'})dW_{t^{'}}+\int_{0}^{t}\frac{b^{*}-b^{2}(t^{'})}{2}dt^{'}}\geq c ~for~ all~ t\right\}>0.
\end{eqnarray*}
We finish the proof.

\subsection{Proof of Theorem \ref{blow up initial value 2}}
The proof of Theorem  \ref{blow up initial value 2} is similar to that of Theorem  \ref{blow up initial value}. We first prove the following lemma.
\begin{lem}\label{p7.2}
Let $s>5/2$ and $b(t)$ satisfy Assumption \ref{Linear noise}. Assume $(u_{0}, \gamma_{0})$ is an $H^{s}\times H^{s}$-valued $\mathcal{F}_{0}$-measurable random variable. Let $K=\frac{\sqrt{2}}{2} (\|u_{0}\|_{H^{1}}^{2}+\|\gamma_{0}\|_{H^{1}}^{2})^{\frac{1}{2}}$. Then for $v_{1}, v_{2}$ defined by \eqref{v1,v2 4}, \eqref{v1,v2 5},
\begin{eqnarray*}
N(\omega,t):=\int_{\mathbb R}v_{1x}^{3}(\omega,t, q(\omega, t,x))dx
\end{eqnarray*}
satisfies $\mathbb P-a.s.$

\begin{eqnarray}\label{7.19}
\frac{d}{dt}N(t)\leq \frac{15\beta}{4}K^{4}-\frac{\beta}{4K^{2}}N^{2}(t),
 \ \ t<\tau^\ast.
\end{eqnarray}
Moreover, if $\mathbb P-a.s.$ $N(0)<-\sqrt{15}K^{3}$, then $\mathbb P-a.s.$. $N(t)$ is non-increasing on $[0, \tau^{\ast})$ and
\begin{eqnarray}\label{7.20}
N(t)<-\sqrt{15}K^{3}, \ \ t\in[0, \tau^{\ast}).
\end{eqnarray}
\begin{pf}
Differentiating the first equation of \eqref{v1,v2 2} with respect to $x$, and using the $\partial_{x}^{2}(1-\partial_{x}^{2})^{-1}f=\partial_{x}^{2}G*f=G*f-f$, we have
\begin{eqnarray}\label{3.2}
v_{1xt}+\frac{\beta}{2}v_{1x}^{2}+\beta v_{1}v_{1xx}+\beta G\ast\left(v_{1}^{2}+\frac{1}{2}v_{1x}^{2}+\frac{1}{2}v_{2}^{2}-\frac{1}{2}v_{2x}^{2}\right)-\beta\left( v_{1}^{2}+\frac{1}{2}v_{2}^{2}-\frac{1}{2}v_{2x}^{2}\right)=0.
\end{eqnarray}
Let $N(t):=\int_{\mathbb R}v_{1x}^{3}(\omega, t, x)dx, t\geq0$. Multiplying \eqref{3.2} with $v_{1x}^{2}$ and integrating by parts subsequently, by $G*(v_{1}^{2}+\frac{1}{2}v_{1x}^{2})(x)\geq \frac{1}{2}v_{1}^{2}(x)$, we get
\begin{align*}
\frac{1}{3}\frac{dN(t)}{dt}
=&-\frac{\beta}{6}\int_{\mathbb R}v_{1x}^{4}dx-\beta \int_{\mathbb R}v_{1x}^{2} G\ast(v_{1}^{2}+\frac{1}{2}v_{1x}^{2}+\frac{1}{2}v_{2}^{2}-\frac{1}{2}v_{2x}^{2})dx+\beta\int_{\mathbb R}v_{1x}^{2}(v_{1}^{2}+\frac{1}{2}v_{2}^{2}-\frac{1}{2}v_{2x}^{2})dx\nonumber\\
\leq&-\frac{\beta}{6}\int_{\mathbb R}v_{1x}^{4}dx+\frac{\beta}{2}\int_{\mathbb R}v_{1}^{2}v_{1x}^{2}dx+\frac{\beta}{2}\int_{\mathbb R}v_{1x}^{2}G*v_{2x}^{2}dx+\frac{\beta}{2}\int_{\mathbb R}v_{1x}^{2}v_{2}^{2}dx\nonumber\\
\leq&-\frac{\beta}{6}\int_{\mathbb R}v_{1x}^{4}dx+\frac{\beta}{2}\int_{\mathbb R}v_{1x}^{2}(v_{1}^{2}+v_{2}^{2})dx+\frac{\beta}{4}\|v_{2x}^{2}\|_{L^{1}}\int_{\mathbb R}v_{1x}^{2}dx.
\end{align*}
In view of Sobolev's embedding and the invariant property of $\|v_{1}(t)\|^{2}_{H^{1}}+\|v_{2}(t)\|^{2}_{H^{1}}=\|u_{0}\|^{2}_{H^{1}}+\|\gamma_{0}\|^{2}_{H^{1}}$, we find that
\begin{eqnarray*}
\frac{3}{2}\int_{\mathbb R}v_{1x}^{2}(v_{1}^{2}+v_{2}^{2})dx+\frac{3}{4}\|v_{2x}^{2}\|_{L^{1}}\int_{\mathbb R}v_{1x}^{2}dx\leq\frac{3}{4}(\|u_{0}\|^{2}_{H^{1}}+\|\gamma_{0}\|^{2}_{H^{1}})^{2}+\frac{3}{16}
(\|u_{0}\|^{2}_{H^{1}}+\|\gamma_{0}\|^{2}_{H^{1}})^{2}.
\end{eqnarray*}
On the other hand, the Cauchy-Schwarz inequality implies that
\begin{eqnarray*}
\bigg|\int_{\mathbb R}v_{1x}^{3}dx\bigg|\leq\bigg(\int_{\mathbb R}v_{1x}^{4}dx\bigg)^{\frac{1}{2}}\bigg(\int_{\mathbb R}v_{1x}^{2}dx\bigg)^{\frac{1}{2}},
\end{eqnarray*}
hence,
\begin{eqnarray*}
\int_{\mathbb R}v_{1x}^{4}dx\geq\frac{1}{\|u_{0}\|^{2}_{H^{1}}+\|\gamma_{0}\|^{2}_{H^{1}}}\bigg(\int_{\mathbb R}v_{1x}^{3}dx\bigg)^{2}.
\end{eqnarray*}
As defined in Lemma 6.5, $K=\frac{\sqrt{2}}{2} (\|u_{0}\|_{H^{1}}^{2}+\|\gamma_{0}\|_{H^{1}}^{2})^{\frac{1}{2}}$, we have the similar Riccati type equation
\begin{eqnarray*}
\frac{d N(t)}{dt}\leq-\frac{\beta}{4K^{2}}N^{2}(t)+\frac{15\beta}{4}K^{4},
\end{eqnarray*}
which is \eqref{7.19}. In order to prove \eqref{7.20}, define stopping time
\begin{eqnarray*}
\chi(w):=\inf\bigg\{t\in[0,\tau^\ast): N(w, t)>-\sqrt{15}K^{3}\bigg\}.
\end{eqnarray*}
If $N(0)<-\sqrt{15}K^{3}$, then $\mathbb P\{\chi(\omega)>0\}=1$. From the definition of $\chi(\omega)$, we find that $w\in\Omega$, $\chi(\omega)\leq\tau^{*}$. From \eqref{7.19}, we conclude that $N(\omega, t)$ is nonincreasing for $t\in[0, \chi(\omega))$. Hence by the continuity of the path of $N(\omega, t)$, we obtain that $N(\omega, t)\leq N(0)<-\sqrt{15}K^{3}$. In view of the time continuity of $N(\omega, t)$ again, we find that $\mathbb P\{\chi=\tau^{*}\}=1.$ Therefore, \eqref{7.20} is true.
\end{pf}
\end{lem}
\textbf{Proof of Theorem \ref{blow up initial value 2}.} From Lemma \ref{p7.2} and \eqref{u0x1}, we rewrite \eqref{7.19} as
\begin{eqnarray*}
\frac{d}{dt}N(t)&\leq& -\frac{\beta(t)}{4K^{2}}\bigg(1-\frac{15K^{6}}{N^{2}(0)}\bigg)N^{2}(t)-\left(\frac{N^{2}(t)}{N^{2}(0)}-1\right)\frac{15\beta(t)}{4}K^{4}\nonumber\\
&\leq& -\frac{\beta(t)}{4K^{2}}\bigg(1-\frac{15K^{6}}{N^{2}(0)}\bigg)N^{2}(t),\ \ t\in[0, \tau^{\ast}).
\end{eqnarray*}
Integrating on both sides leads to $\mathbb P-a.s.$
\begin{eqnarray*}
\frac{1}{N(t)}-\frac{1}{N(0)}\geq \left(1-\frac{15K^{6}}{N^{2}(0)}\right)\int_{0}^{t}\frac{\beta(t^{'})}{4K^{2}}dt^{'},~~t\in[0, \tau^{*}).
\end{eqnarray*}
Assuming $\Omega^{'}=\{\omega: \beta(t,\omega)\geq c e^{-\frac{b^{*}}{2}t}~for~ all~ t\}$, also due to $N(t)<-\sqrt{15}K^{3}$, we get $\mathbb{P}-a.s.$ $\omega\in\Omega^{'}$
\begin{eqnarray*}
-\frac{1}{N(0)}\geq \bigg(\frac{1}{4K^{2}}-\frac{15K^{4}}{4N^{2}(0)}\bigg)\int_{0}^{\tau^{*}}\beta(t^{'})dt^{'}\geq\bigg(\frac{1}{4K^{2}}-\frac{15K^{4}}{4N^{2}(0)}\bigg)\bigg(\frac{2c}{b^{*}}-\frac{2c}{b^{*}}e^{-\frac{b^{*}}{2}\tau^{*}}\bigg).
\end{eqnarray*}
If  $N(0)<-\sqrt{\frac{(b^{*})^{2}K^{4}}{c^{2}}+15 K^{6}}-\frac{b^{*}K^{2}}{c}$, we obtain on $\Omega'$
\begin{eqnarray*}
\bigg(\frac{1}{4K^{2}}-\frac{15K^{4}}{4N^{2}(0)}\bigg)\frac{2c}{b^{*}}e^{-\frac{b^{*}}{2}\tau^{*}} \geq\frac{2c}{b^{*}}\bigg(\frac{1}{4K^{2}}-\frac{15K^{4}}{4N^{2}(0)}\bigg)+\frac{1}{N(0)}>0.
\end{eqnarray*}
Therefore we obtain $\tau^{*}<\infty$ $\mathbb{P}-a.s.$ on $\Omega^{'}$, which means that
\begin{eqnarray*}
\mathbb P\{\tau^{*}<\infty\}\geq\mathbb P\{\beta(t)\geq c e^{-\frac{b^{*}}{2}t}~for~ all~ t\}=\mathbb P\left\{e^{\int_{0}^{t}b(t^{'})dW_{t^{'}}+\int_{0}^{t}\frac{b^{*}-b^{2}(t^{'})}{2}dt^{'}}\geq c ~for~ all~ t\right\}>0.
\end{eqnarray*}
So, the proof is finished.

\section{Acknowledgments}
This paper is supported by Fundamental Research Funds for the Central Universities (No. 22D110913).

\end{document}